\documentclass[10pt,leqno]{article}
\usepackage[OT1]{fontenc}
\usepackage[utf8]{inputenc}
\usepackage[
  paperwidth=486bp,
  paperheight=720bp,
  textwidth=369bp,
  textheight=600bp,
  top=60bp,
  hcentering,
  footskip=24bp
]{geometry}
\usepackage{amssymb,mathtools,bm}
\usepackage{amsthm}
\newtheoremstyle{boldtheorem}
  {\topsep}{\topsep}
  {\itshape}{}
  {\bfseries}{.}{5pt}
  {\thmname{#1}\thmnumber{ #2}\thmnote{ \textbf{(#3)}}}
\theoremstyle{boldtheorem}
\newtheorem{lemma}{Lemma}[section]
\newtheorem{theorem}{Theorem}[section]

\theoremstyle{remark}
\newtheorem{remark}{Remark}[section]
\makeatletter
\renewenvironment{proof}[1][\proofname]{%
  \par\pushQED{\qed}%
  \normalfont\topsep6\p@\@plus6\p@\relax
  \trivlist
  \item[\hskip\labelsep\normalfont\bfseries #1\@addpunct{:}]\ignorespaces
}{%
  \popQED\endtrivlist\@endpefalse
}
\makeatother
\numberwithin{equation}{section}
\usepackage{graphicx}
\usepackage{booktabs,tabularx,array,multirow,adjustbox}
\usepackage{subcaption}
\usepackage{float}
\usepackage{placeins,needspace}
\usepackage{hyperref}
\hypersetup{hidelinks}
\usepackage{microtype}
\newcolumntype{Y}{>{\raggedright\arraybackslash}X}
\floatstyle{ruled}
\ifdefined\algorithm\else\newfloat{algorithm}{tbp}{loa}\fi
\floatname{algorithm}{Algorithm}
\ifdefined\procedure\else\newfloat{procedure}{tbp}{lop}\fi
\floatname{procedure}{Procedure}
\usepackage{algpseudocode}
\providecommand{\algcall}[1]{\operatorname{\text{#1}}}

\usepackage{titlesec}
\titleformat{\section}[runin]
  {\normalfont\normalsize\bfseries}{\thesection.}{0.5em}{}[.]
\titlespacing*{\section}{\parindent}{1.5ex plus .5ex minus .2ex}{0.5em}
\titleformat{\subsection}[runin]
  {\normalfont\normalsize\bfseries}{\thesubsection.}{0.5em}{}[.]
\titlespacing*{\subsection}{\parindent}{1.2ex plus .5ex minus .2ex}{0.5em}
\titleformat{\subsubsection}[runin]
  {\normalfont\normalsize\bfseries}{\thesubsubsection.}{0.5em}{}[.]
\titlespacing*{\subsubsection}{\parindent}{1ex plus .5ex minus .2ex}{0.5em}
\DeclareCaptionFormat{papertable}{%
  {\centering #1\par}\nobreak
  \vspace{2pt}%
  \hspace*{18pt}#3\par}
\DeclareCaptionFormat{paperfigure}{%
  \hspace*{18pt}#1#2#3\par}
\usepackage{etoolbox}
\AtBeginEnvironment{algorithmic}{\footnotesize\setlength{\baselineskip}{10pt}}

\makeatletter
\renewcommand{\@maketitle}{%
  \newpage
  \null
  \vskip 0.25em
  \begin{center}
    {\fontsize{12}{14}\selectfont\bfseries
     \MakeUppercase{\@title}\par}
    \vskip 11pt
    {\normalsize\MakeUppercase{\@author}\par}
  \end{center}
  \vspace{0.2em}
}
\makeatother
\renewenvironment{abstract}{%
  \par\noindent\hspace*{\parindent}%
  \textbf{Abstract.}\enspace\ignorespaces
}{%
  \par\vspace{1.0\baselineskip}
}
\title{A Hyper-Reduced Reduced Order Model Built From Streaming Tensor Train Data}

\author{Xiaoshan Wu\thanks{\NoCaseChange{School of Mathematics, Harbin Institute of Technology, Harbin, China. Email: \href{mailto:zoew5622@gmail.com}{zoew5622@gmail.com}.}}, Zhichao Peng\thanks{\NoCaseChange{Department of Mathematics, The Hong Kong University of Science and Technology, Clear Water Bay, Kowloon, Hong Kong, China. Email: \href{mailto:pengzhic@ust.hk}{pengzhic@ust.hk}.}}}
\date{November 2026}

\begin{document}
\maketitle

\begin{abstract}
Low-rank PDE solvers based on tensor-decompositions inside physical space leverages low-rank structures across different physical dimensions, while projection-based reduced order models (ROMs) for parametric PDEs exploit low-rank structures across parameter spaces. Recently, to use the complimentary strength of these two methods, a ROM framework  is proposed by incrementally constructing compressed reduced bases directly from streaming tensor-train (TT) snapshots provided by low-rank solvers, but it omits hyper reduction and consequently incur substantial online costs for nonlinear or non-affine parametric PDEs. To address this issue, we develop a hyper-reduced TT ROM for streaming TT data by combining incremental low-rank TT compression and nested TT discrete empirical interpolation (TT-DEIM).  Numerical experiments on nonlinear and non-affine parametric PDEs achieve mean relative errors on the order of $10^{-5}$–$10^{-4}$ relative to full-order reference solutions, matching the accuracy of non-hyper-reduced TT-Galerkin ROMs while achieving mean online speedups of approximately $8\times$ to $1500\times$ comapred to them.
\end{abstract}

\section{Introduction}
\label{sec:1}

Multi-query applications, including uncertainty quantification, inverse problems, and design optimization, often require repeated solutions of a large-scale parameteric partial differential equation (PDE) at different parameter values, incurring substantial computational costs, espiecially for high deimnsional PDEs. Projection-based reduced order models (ROMs) \cite{benner2015survey} and low-rank solvers based on tensor or matrix decompositions \cite{bachmayr2023lowrank,einkemmer2025review} are two dimensionality reduction techniques to reduce the costs for such problems.

Tensor-based low-rank solvers exploit low-rank structure across different dimensions of the PDE solution through matrix or tensor compressions~\cite{bachmayr2023lowrank}. Notable applications include kinetic equations~\cite{einkemmer2025review}, quantum
dynamics~\cite{meyer1990multiconfigurational,wang2003multilayer}, optimal control~\cite{dolgov2021hjb}, and parametric PDEs with random coefficients in uncertainty quantification~\cite{dolgov2019hybrid}. For time integration, dynamical low-rank approximation (DLRA) evolves the solution  along the tangent space of the low-rank manifold~\cite{koch2007dynamical,lubich2013dynamical}, whereas step-and-truncate (SAT) schemes compress the solution after each discrete time step~\cite{kormann2015semilagrangian,rodgers2022adaptive,guo2022lowrank}.  Cross tenosr approximation, reconstructing low-rank representations from selected entries~\cite{oseledets2010ttcross}, has also been leveraged in recent sampling-based low-rank PDE solvers ~\cite{ghahremani2024cross,dektor2025collocation,zheng2025slar}.

Projection-based ROMs exploits low-rank structures across parameters of parametric PDEs~\cite{benner2015survey}. It typically follows an offline-online decomposition. Its offline stage constructs low-dimensional reduced space by extracting low-rank structures from solution snapshots. Its online stage predicts low-fidelity solutions for new parameters efficiently through projections onto the reduced space constructed offline. Representative methods include proper orthogonal decomposition (POD)~\cite{sirovich1987pod,kunisch2001galerkin}, reduced basis methods \cite{buffa2012greedy,haasdonk2008reduced}, balance truncation \cite{moore1981principal,glover1984all}, and $H_2$ rational interpolation \cite{gugercin2008h2}.

However, when a PDE contains non-affine parameter dependencies or nonlinear terms, projection alone no longer ensures online efficiency because assembling the reduced operators still requires full-order evaluations of non-affine or nonlinear terms. To restore online complexity independent of the full-order dimension,  hyper reduction approximates the non-affine or nonlinear terms through
evaluations at a small set of selected entries \cite{ryckelynck2005hyperreduction,farhat2020computational}.  Interpolation-based approaches, such as the empirical interpolation method (EIM)~\cite{barrault2004eim} and its discrete counterpart (DEIM)~\cite{chaturantabut2010deim} approximate these terms in an empirical basis, and the corresponding reduced operators are assembled online via interpolations at selected spatial indices~\cite{grepl2007nonaffine}. Alternative hyper reduction strategies include Gauss--Newton with approximated tensors (GNAT)~\cite{carlberg2013gnat}, which samples residual and Jacobian rows, energy-conserving sampling and weighting (ECSW)~\cite{farhat2015ecsw}, which sums weighted contributions over a subset of mesh elements, and the LP empirical quadrature method~\cite{yano2019lp}, which uses an offline linear program to select points and positive weights for evaluating reduced residual integrals online.

Tensor decompositions have also been used to build ROMs from multidimensional snapshot data. Tensorial ROMs compress snapshots across time and parameters to construct parameter-specific reduced spaces~\cite{mamonov2022interpolatory,mamonov2024tensorial}, while TT-ROM applies the tensor-train (TT) decomposition~\cite{oseledets2011tt} to build bases for the solution, residual, and Jacobian on Cartesian grids~\cite{mueller2026ttrom}.  However, both tensorial ROM and TT-ROM begin by assembling full-order snapshot tensors; for high-dimensional kinetic equations and many-body quantum dynamics, storing such high-dimensional tensors becomes memory-expensive.
To address this limitation, the low-rank incremental TT compression is developed \cite{shang2026incremental} to directly build a compressed reduced basis for a parametric kinetic equation from streaming low-rank data, exploiting the memory and computational benefits of tensor-based low-rank solvers. However, the resulting ROM lacks low-rank-compatible hyper reduction and remains restricted to linear PDEs with affine parameter dependence.

We extend the streaming TT-ROM~\cite{shang2026incremental} to parametric PDEs with non-affine  or nonlinear terms by combining incremental TT reduced basis construction with nested TT-DEIM~\cite{dektor2025collocation}. Offline, we incrementally build the solution and solution basis from the TT-format outputs of low-rank solvers. Nested TT-DEIM selects a fixed set of interpolation points directly from the cores of the empirical basis. Online, at a new parameter value, we evaluate non-affine or nonlinear terms at these selected points and cheaply constructions reduced operators through interpolations.
The resulting hyper-reduced TT ROM preserves the memory advantage of low-rank solvers during offline training by building both the solution and empirical bases without assembling full-tensor snapshots.
Through a series of numerical tests, we demonstrate that hyper reduction in TT-ROMs preserves the order of prediction error and yields substantial online acceleration over the TT-ROMs without it, achieving speedups exceeding $1500\times$ in a non-affine radiative transfer benchmark.

The rest of the manuscript is organized as follows.
Section~\ref{sec:2} briefly reviews low-rank sovlers, projection-based ROMs, hyper reductions and recent TT-based ROMs.  Section~\ref{sec:3} presents the hyper-reduced TT-ROM, and estimates the computational cost and storage of its offline stage. Section~\ref{sec:4} bounds the TT-DEIM interpolation error on the training set. Section~\ref{sec:5} reports numerical results, and Section~\ref{sec:6} concludes the paper.

\section{Research Background and Motivation}
\label{sec:2}

We first consider parametric PDEs posed in a $d$-dimensional physical space arising from multiquery applications such as sensitivity analysis, uncertainty quantification or optimization.  After discretization, suppose there is $n$ DOFs in each dimension. Then, the dsicrete solution solution can be represented by \(\boldsymbol{u}(\boldsymbol{\mu})\in\mathbb R^N\), where \(N=n^d\) is the number of full-order degrees of freedom (DOFs) and $\boldsymbol{\mu}$ denotes the collection of parameters of the underlying parametric PDE such as material or boundary configurations. We write the possibly nonlinear discrete problem in residual form as
\begin{equation}
\label{eq:s2-d1}
\boldsymbol{r}(\boldsymbol{u}(\boldsymbol{\mu});\boldsymbol{\mu})=\boldsymbol{0},\qquad
\boldsymbol{r}(\cdot;\boldsymbol{\mu}):\mathbb R^N\to\mathbb R^N.
\end{equation}
For time-dependent problems, this formulation can be applied at each implicit time step treating previous solutions as known sources.

When the dimension $d$ is high, solving the PDE for even one parameter can lead to substantial computational and storage costs. To address the computational burden in multiquery applications, two 
approaches have been developed. The first is data-driven reduced order models (ROMs) ~\cite{benner2015survey} that builds low-dimensional linear approximations of the PDE solution by extracting low-rank structures across parameters from solution data. The second is low-rank numerical solvers~\cite{bachmayr2023lowrank} based on matrix or tensor compressions of the PDE solutions, leveraging the low-rank structures across different physical dimensions. In this section, we first review the basic ideas of low-rank numerical solvers, data-driven ROMs, and recent low-rank ROMs built from streaming low-rank data provided by low-rank solvers, then point out the bottleneck of the existing approach for parametric PDEs with non-affine and nonlinear terms.

\subsection{Low-Rank Numerical Solvers}
\label{sec:2-1}

At a fixed parameter value \(\boldsymbol{\mu}\), low-rank numerical solvers exploit correlations across physical dimensions to compute PDE solutions directly in compressed matrix or tensor formats~\cite{bachmayr2023lowrank}. For example, low-rank solvers for kinetic equations~\cite{einkemmer2021transport,guo2026sweep,han2026adaptive} separate the spatial and velocity variables by approximating the discrete particle distribution
\[
\boldsymbol{F}(\boldsymbol{\mu})
=
\begin{bmatrix}
\boldsymbol{f}(\cdot,v_1;\boldsymbol{\mu}) & \cdots & \boldsymbol{f}(\cdot,v_{N_v};\boldsymbol{\mu})
\end{bmatrix}
\in\mathbb R^{N_x\times N_v}
\]
using a low-rank matrix factorization. Suppressing the dependence on \(\boldsymbol{\mu}\), we write
\begin{equation}
\label{eq:s2-matrix-lowrank}
\boldsymbol{F}\approx \boldsymbol{X}\boldsymbol{K}\boldsymbol{V}^T,
\qquad
\boldsymbol{X}\in\mathbb R^{N_x\times r},\quad
\boldsymbol{K}\in\mathbb R^{r\times r},\quad
\boldsymbol{V}\in\mathbb R^{N_v\times r},
\end{equation}
where \(N_x\) and \(N_v\) denote the numbers of spatial and velocity degrees of freedom, respectively, and \(r\ll\min(N_x,N_v)\). The columns of \(\boldsymbol{X}\) and \(\boldsymbol{V}\) form spatial and velocity bases, while the core matrix \(\boldsymbol{K}\) contains their coupling coefficients. Alternatively, low-rank tensor compressions, such as tensor train (TT)~\cite{oseledets2011tt}, Tucker~\cite{tucker1966notes}, and hierarchical Tucker~\cite{hackbusch2009scheme}, can also be leveraged in PDE solves, for example in \cite{oseledets2012linear,dolgov2014alternating,wang2026dynamical}. When the discretized operators admit compatible low-rank tensor representations, their application can be implemented through tensor contractions. For nonlinear terms challenging to direct tensor decompositions, cross approximation~\cite{oseledets2010ttcross} can be used to construct such representations by evaluating selected entries, and used in PDE solves~\cite{dolgov2019hybrid,ghahremani2024cross,zheng2025slar,dektor2025collocation,appelo2026lraa}.

\subsection{Standard ROM and Hyper Reduction}
\label{sec:2-2}
A complementary dimensionality reduction technique to  is the reduced order models (ROMs)~\cite{benner2015survey}. Here, we briefly review standard proper orthogonal decomposition (POD) following an offline-online decomposition framework. 

\textbf{Offline stage.} In the offline stage, POD constructs a common reduced basis by extracting low-rank structure across parameter values from solution snapshots. Let \(\mathcal{P}\) denote the parameter domain and \(\{\boldsymbol{\mu}_j\}_{j=1}^M\subset\mathcal{P}\) a training set. The corresponding solutions are collected as columns of the snapshot matrix
\begin{equation}
\label{eq:s2-d2}
\boldsymbol{S}=[\boldsymbol{u}(\boldsymbol{\mu}_1),\ldots,\boldsymbol{u}(\boldsymbol{\mu}_M)]\in\mathbb R^{N\times M}.
\end{equation}
Given the thin SVD, \(\boldsymbol{S}=\boldsymbol{L}\boldsymbol{\Sigma} \boldsymbol{R}^T\), with singular values \(\sigma_1\geq\cdots\geq\sigma_\ell\geq0\) and \(\ell=\min(N,M)\), POD retains the leading \(r_u\) left singular vectors as the columns of \(\boldsymbol{V}\). This gives the basis--coefficient representation
\begin{equation}
\label{eq:s2-pod-factorization}
\boldsymbol{S}\approx \boldsymbol{V}\boldsymbol{C},\qquad \boldsymbol{V}^T\boldsymbol{V}=\boldsymbol{I}_{r_u},\qquad
\boldsymbol{C}=\boldsymbol{V}^T\boldsymbol{S}=\boldsymbol{\Sigma}_{1:r_u,1:r_u}\boldsymbol{R}_{:,1:r_u}^T.
\end{equation}
Here \(\boldsymbol{C}\in\mathbb R^{r_u\times M}\) stores the coefficients of the training snapshots in the reduced basis. The dominant singular directions minimize the projection error over all \(r_u\)-dimensional subspaces~\cite{sirovich1987pod,benner2015survey}:
\begin{equation}
\label{eq:s2-pod-error}
\begin{aligned}
\|\boldsymbol{S}-\boldsymbol{V}\boldsymbol{V}^T\boldsymbol{S}\|_F^2
&=\min_{\substack{\boldsymbol{Z}\in\mathbb R^{N\times r_u}\\\boldsymbol{Z}^T\boldsymbol{Z}=\boldsymbol{I}_{r_u}}}
\|\boldsymbol{S}-\boldsymbol{Z}\boldsymbol{Z}^T\boldsymbol{S}\|_F^2=\sum_{k=r_u+1}^{\ell}\sigma_k^2.
\end{aligned}
\end{equation}
For a prescribed relative tolerance \(\varepsilon_{\mathrm{POD}}\), one decide the reduced space dimension by finding the smallest \(r_u\) satisfying \(\sum_{k=r_u+1}^{\ell}\sigma_k^2\leq\varepsilon_{\mathrm{POD}}^2\|\boldsymbol{S}\|_F^2\). 

\textbf{Online stage.}  In the online stage, solutions for a new parameter can be predicted through Galerkin or Petrov-Galerkin projections of the full-order discrete system onto the reduced space constructed offline.
Galerkin projection finds a reduced solution by 
\begin{equation}
\label{eq:s2-d8}
\boldsymbol{u}(\boldsymbol{\mu})\approx \boldsymbol{V}\boldsymbol{c}_u(\boldsymbol{\mu}),\quad\text{where}\; \boldsymbol{V}^T\boldsymbol{r}(\boldsymbol{V}\boldsymbol{c}_u(\boldsymbol{\mu});\boldsymbol{\mu})=\boldsymbol{0}.
\end{equation}
Here \(\boldsymbol{V}\in\mathbb R^{N\times r_u}\) is the reduced basis, and \(\boldsymbol{c}_u(\boldsymbol{\mu})\in\mathbb R^{r_u}\) is the reduced coefficient.

\textbf{Affine offline--online decomposition.} For a linear problem, $\boldsymbol{R}(\boldsymbol{u};\boldsymbol{\mu})=\boldsymbol{A}(\boldsymbol{\mu})\boldsymbol{u}-\boldsymbol{f}(\boldsymbol{\mu})$, with $\boldsymbol{A}(\boldsymbol{\mu})\in\mathbb R^{N\times N}$ and $\boldsymbol{f}(\boldsymbol{\mu})\in\mathbb R^N$, and the reduced equation becomes $\boldsymbol{V}^T\boldsymbol{A}(\boldsymbol{\mu})\boldsymbol{V}\boldsymbol{c}_u(\boldsymbol{\mu})=\boldsymbol{V}^T\boldsymbol{f}(\boldsymbol{\mu})$. Suppose further that the operator and right-hand side admit affine decompositions with respect to the parameter $\boldsymbol{\mu}$:
\begin{equation}
\label{eq:s2-d9}
\boldsymbol{A}(\boldsymbol{\mu})=\sum_{q=1}^{Q_A}\theta_q^A(\boldsymbol{\mu})\boldsymbol{A}_q,\quad
\boldsymbol{f}(\boldsymbol{\mu})=\sum_{q=1}^{Q_f}\theta_q^f(\boldsymbol{\mu})\boldsymbol{f}_q,
\end{equation}
where \(\boldsymbol{A}_q\) and \(\boldsymbol{f}_q\) are parameter-independent. The reduced operator and projected right-hand side can then be fast assembled online as linear combinations of precomputed projections offline:
\begin{equation}
\boldsymbol{A}_{r,q}=\boldsymbol{V}^T\boldsymbol{A}_q\boldsymbol{V},\qquad \boldsymbol{f}_{r,q}=\boldsymbol{V}^T\boldsymbol{f}_q.
\end{equation}

\textbf{Hyper reduction.} A general residual $\boldsymbol{r}(\boldsymbol{u};\boldsymbol{\mu})$ may involve  nonlinear or non-affine terms. Projecting such terms can still require full-order evaluations, offsetting the saving of approximating the high-dimensional full-order system with a low-dimensional reduced-order system. Hyper reduction techniques address this bottleneck by approximating these terms using selected interpolation points ~\cite{barrault2004eim,chaturantabut2010deim}.

One popular approach is the discrete empirical interpolation method (DEIM), which constructs a reduced basis (the DEIM basis) for nonlinear or non-affine terms and selects interpolation indices for fast online evaluation. To illustrate, let $\boldsymbol{g}(\boldsymbol{u};\boldsymbol{\mu})\in\mathbb R^N$ be a vector term in the residual. Offline, we form the snapshot matrix $\boldsymbol{S}_g=[\boldsymbol{g}^{(1)},\ldots,\boldsymbol{g}^{(M)}]$, with $\boldsymbol{g}^{(j)}=\boldsymbol{g}(\boldsymbol{u}(\boldsymbol{\mu}_j);\boldsymbol{\mu}_j)$, and take its leading $m$ left singular vectors, $\boldsymbol{W}=[\boldsymbol{w}_1,\ldots,\boldsymbol{w}_m]\in\mathbb R^{N\times m}$, as the DEIM basis. The interpolation points are represented by a sampling matrix $\boldsymbol{P}=[\boldsymbol{e}_{p_1},\ldots,\boldsymbol{e}_{p_m}]\in\mathbb R^{N\times m}$ (where $\boldsymbol{e}_i$ denotes the $i$th canonical unit vector), which is constructed greedily: at each step $j$, $p_j$ selects the largest-magnitude entry of the residual obtained by interpolating $\boldsymbol{w}_j$ using the preceding basis vectors and indices.
Online, DEIM approximates \(\boldsymbol{g}(\boldsymbol{V}\boldsymbol{c}_u;\boldsymbol{\mu})\) using its entries at the selected indices \(\mathcal{I}=\{p_1,\ldots,p_m\}\). Suppressing the state and parameter arguments and using MATLAB notation, DEIM yields \begin{equation}
\label{eq:s2-deim-preprojection}
\begin{aligned}
\boldsymbol{V}^T\widehat{\boldsymbol{g}}
=(\boldsymbol{V}^T\boldsymbol{W})(\boldsymbol{P}^T\boldsymbol{W})^{-1}\boldsymbol{P}^T\boldsymbol{g}=(\boldsymbol{V}^T\boldsymbol{W})[\boldsymbol{W}(\mathcal{I},:)]^{-1}\boldsymbol{g}(\mathcal{I}).
\end{aligned}
\end{equation}
Precomputing \(\boldsymbol{V}^T\boldsymbol{W}\) and a factorization of \(\boldsymbol{P}^T\boldsymbol{W}\) offline reduces online evaluation to evaluating $\boldsymbol{g}$ at selected indices $\mathcal{I}$ and performing low-dimensional  operations.  Algorithm~\ref{alg:deim} summarizes the DEIM algorithm following ~\cite{chaturantabut2010deim}.

\begin{algorithm}[htbp]
\caption{DEIM interpolation index selection\cite{chaturantabut2010deim}}
\label{alg:deim}
\begin{algorithmic}[1]
\Require Basis $\boldsymbol{W}=[\boldsymbol{w}_1,\ldots,\boldsymbol{w}_m]\in\mathbb R^{N\times m}$ with linearly independent columns
\Ensure Ordered indices $\mathcal{I}=(p_1,\ldots,p_m)$ and sampling matrix $\boldsymbol{P}$
\State Choose $p_1\in\operatorname*{arg\,max}_{1\leq i\leq N}|\boldsymbol{W}(i,1)|$; set $\mathcal{I}\gets(p_1)$
\For{$j=2,\ldots,m$}
  \State Solve $\boldsymbol{W}(\mathcal{I},1:j-1)\boldsymbol{c}=\boldsymbol{W}(\mathcal{I},j)$
  \State $\boldsymbol{z}\gets \boldsymbol{W}(:,j)-\boldsymbol{W}(:,1:j-1)\boldsymbol{c}$
  \State Choose $p_j\in\operatorname*{arg\,max}_{i\notin\mathcal I}|z_i|$
  \State Append $p_j$ to $\mathcal{I}$
\EndFor
\State $\boldsymbol{P}\gets[\boldsymbol{e}_{p_1},\ldots,\boldsymbol{e}_{p_m}]$
\State \Return $\mathcal{I},\boldsymbol{P}$
\end{algorithmic}
\end{algorithm}

\subsection{TT-ROM from Full-Order Snapshot Data}
\label{sec:2-3}

In \cite{mueller2026ttrom}, a TT-ROM is developed from full-order snapshot tensors, exploiting Cartesian tensor
structure in basis construction and hyper reduction. 

The TT format approximates \(\mathcal{X}\in\mathbb R^{n_1\times\cdots\times n_d}\) by contracting three-dimensional cores~\cite{oseledets2011tt}:
\begin{equation}
\label{eq:s2-tt-representation}
\begin{gathered}
\mathcal{X}(i_1,\ldots,i_d)\approx
\mathcal{G}_1(1,i_1,:)\mathcal{G}_2(:,i_2,:)\cdots \mathcal{G}_d(:,i_d,1),\\
\mathcal{G}_k\in\mathbb R^{r_{k-1}\times n_k\times r_k},\qquad r_0=r_d=1.
\end{gathered}
\end{equation}
The TT ranks \(r_k\) determine the storage cost
\(\sum_{k=1}^d n_k r_{k-1}r_k\), bounded by
\(\mathcal{O}(dn\bar r^2)\) when \(n_k\leq n\) and \(r_k\leq\bar r\).

For basis construction, the full-order solution snapshots are reshaped into an order-\(d\) tensor, and TT-SVD replaces the SVD step in POD to construct a reduced basis in TT format. For hyper reduction, TT-SVD similarly constructs DEIM bases in TT format for residuals and Jacobians, while a variant of TT-DEIM \cite{dektor2025collocation}, namely TT-MDEIM, selects interpolation indices through operations only on TT cores.

\subsection{TT-ROM from Low-Rank Solution Data and Its Limitations}
\label{sec:2-4}
Despite recent advances in TT-ROM~\cite{mueller2026ttrom}, basis construction still starts from full-order solution snapshots and does not directly exploit the compressed outputs of low-rank solvers. To preserve the storage and computational cost benefits of low-rank solvers offline, \cite{shang2026incremental} incrementally constructs a TT reduced basis from streaming TT-formatted snapshots. This approach avoids both reconstructing full tensors and uses  incoming low-rank snapshots in a streaming manner.
However, only affine problems with tensor-structured operators and right-hand sides are considered in~\cite{shang2026incremental}. To maintain the computational benefits of low-rank solvers for non-affine and nonlinear problems, we incoporate  hyper reduction based on TT-DEIM \cite{dektor2025collocation}, similar to \cite{mueller2026ttrom}.

\begin{remark}
Parametric tensor solvers treat parameters as additional tensor modes and solve jointly over physical and parameter variables, avoiding separate snapshot generation~\cite{kressner2011low,dolgov2014alternating}. However, robust convergence with moderate intermediate ranks can be challenging, particularly when preconditioners must remain effective across the parameter domain and support efficient low-rank operations. Even for a fixed parameter value, substantial rank growth can make a TT-based solver slower than its full-rank counterpart~\cite{gorodetsky2025thermal}. Snapshot-based ROMs allow flexible training sets and specialized solvers with preconditioners tailored to individual parameter values. 
\end{remark}

\section{Hyper-Reduced TT-ROM from Streaming TT-Data}
\label{sec:3}

We construct hyper-reduced models for the parametric residual equation~\eqref{eq:s2-d1} directly from streaming TT solution data and sampled nonlinear or non-affine terms. The main task is to build the solution and DEIM spaces, interpolation indices, and reduced projection operators without full tensor reconstruction. Section~\ref{sec:3-1} states the data assumptions and overall procedure. Section~\ref{sec:3-2} presents the streaming basis construction in \cite{shang2026incremental},  Section~\ref{sec:3-3} uses TT-DEIM \cite{dektor2025collocation} to perform hyper reduction. Section~\ref{sec:3-4} gives the online solve. Section~\ref{sec:3-5} summarizes offline storage and construction costs.

\subsection{Data Assumptions and Overall Offline--Online Procedure}
\label{sec:3-1}

Let $\mathcal{P}_{\mathrm{train}}=\{\boldsymbol{\mu}^{(j)}\}_{j=1}^M$ be the training set, with solution snapshots $\boldsymbol{u}^{(j)}$ supplied in tensor-train (TT) format by a low-rank solver. We incrementally build a compressed reduced basis using the low-rank incremental TT compression (LR-iTT) method of~\cite{shang2026incremental}. This approach assumes that each incoming TT snapshot $\mathcal{X} = \mathcal{G}_1 \cdots \mathcal{G}_d$ is right-orthogonal and incrementally updates a left-orthogonal accumulated snapshot TT $\mathcal{S} = \mathcal{U}_1 \cdots \mathcal{U}_d \mathcal{U}_{d+1}$:
\[
\begin{aligned}
\operatorname{unfold}_R(\mathcal{G}_i)\operatorname{unfold}_R(\mathcal{G}_i)^T &= \boldsymbol{I}, && i=2,\ldots,d,\\
\operatorname{unfold}_L(\mathcal{U}_i)^T\operatorname{unfold}_L(\mathcal{U}_i) &= \boldsymbol{I}, && i=1,\ldots,d.
\end{aligned}
\]
For a tensor core $\mathcal{C} \in \mathbb{R}^{a \times n \times b}$, $\operatorname{unfold}_L(\mathcal{C}) \in \mathbb{R}^{(an) \times b}$ and $\operatorname{unfold}_R(\mathcal{C}) \in \mathbb{R}^{a \times (nb)}$. Under these orthogonalizations, the first core $\mathcal{G}_1$ carries the snapshot norm ($\|\mathcal{G}_1\|_F = \|\mathcal{X}\|_F$). These properties are maintained in compressed format via QR sweeps~\cite{oseledets2011tt}, allowing local residual tests to control the global tensor approximation error. Snapshots in this format can be natively or cheaply produced by several low-rank solvers, such as~\cite{guo2026sweep,han2026adaptive}.

For each nonlinear or non-affine target quantity $\boldsymbol{g}(\boldsymbol{u};\boldsymbol{\mu})$, we assume inexpensive entrywise evaluation via
\[
g_{\mathrm{entry}}(\boldsymbol{u},\boldsymbol{\mu},i) = [\boldsymbol{g}(\boldsymbol{u};\boldsymbol{\mu})]_i.
\]
This entrywise evaluation avoids assembling the full tensor quantity, enabling efficient snapshot construction offline and reduced-system assembly online. 

\textbf{Offline stage: basis construction and hyper reduction.}
The offline stage constructs compressed reduced bases, interpolation points, and performs necessary precomputations to enable fast online predictions :
\begin{enumerate}
\item \textbf{Construct the reduced basis.} Incrementally incorporate the solution snapshots into $\mathcal{S}_u$ and extract an orthonormal TT basis (Section~\ref{sec:3-2}).
\item \textbf{Hyper reduction.} For each nonaffine or nonlinear term $\boldsymbol{g}$, construct TT approximations from selected entries, and incrementally compress them into $\mathcal{S}_g$ to extract a DEIM basis following Section \ref{sec:3-2}. Apply TT-DEIM to the compressed DEIM basis to select interpolation indices $\mathcal{I}_g$ core by core. Then, perform necessary precomputations to enable fast online predictions (Section~\ref{sec:3-3}).
\end{enumerate}

\textbf{Online stage: efficient solution at a new parameter.}
In the online stage, we exploit the compressed reduced basis and selected interpolation points offline to efficiently compute reduced solutions for new parameters. For each new parameter $\boldsymbol{\mu}$:
\begin{enumerate}
\item \textbf{Sample and interpolate.} Determine the interpolation coefficients of each target quantity from its selected entries.
\item \textbf{Assemble the reduced residual.} Use these coefficients and offline precomputations to evaluate the hyper-reduced residual.
\item \textbf{Solve and recover the TT solution.} Solve the reduced equation, repeating the sampled evaluations and reduced assembly as needed during nonlinear iteration, and reconstruct the reduced solution in TT format.
\end{enumerate}

\subsection{Offline Reduced-Basis Construction from Streaming TT Data}
\label{sec:3-2}

In the offline stage, we incrementally construct the compressed reduced basis from streaming TT snapshots $\mathcal{X}^{(1)},\ldots,\mathcal{X}^{(M)}$ of physical mode sizes $q_1,\ldots,q_d$ using the LR-iTT-F algorithm in ~\cite{shang2026incremental}. Storing the sample index as an additional mode, the accumulated snapshot tensor after $k$ updates takes the form
\begin{equation}
\label{eq:s3-d16}
\mathcal{S}^{(k)}=\mathcal{U}_1^{(k)}\cdots \mathcal{U}_d^{(k)}\mathcal{U}_{d+1}^{(k)},\qquad
\mathcal{S}^{(k)}(:,\ldots,:,j)\approx \mathcal{X}^{(j)},\quad j=1,\ldots,k,
\end{equation}
where the physical cores $\mathcal{U}_1^{(k)},\ldots,\mathcal{U}_d^{(k)}$ encode a compressed basis of the physical modes and the terminal core $\mathcal{U}_{d+1}^{(k)}$ stores the corresponding coefficients for each parameter. At each step, the algorithm updates the physical cores by adding only those components of the incoming snapshot that are not well-represented by the existing basis, governed by a prescribed local residual tolerance.

We initialize $\mathcal{S}^{(1)}$ from the first snapshot,  by applying a left-orthogonalization if necessary.  For each subsequent snapshot, assumed to be in the right-orthogonal form specified in Section~\ref{sec:3-1}, we denote its cores by $\mathcal{X}_i\in\mathbb{R}^{s_{i-1}\times q_i\times s_i}$ ($s_0=s_d=1$), set $\boldsymbol{W}_1=\operatorname{unfold}_L(\mathcal{X}_1)$, and sweep across the physical modes $i=1,\ldots,d$, performing the following operations at each mode:
\begin{enumerate}
\item \textbf{Identify new components via orthogonal projection.} Let $\widetilde{\boldsymbol{U}}_i = \operatorname{unfold}_L(\mathcal{U}_i^{(k-1)})$ denote the left unfolding of the $i$-th core of the accumulation tensor, zero-padded along its left rank if the preceding mode's rank increased. Project $\boldsymbol{W}_i$ onto the orthogonal complement of $\widetilde{\boldsymbol{U}}_i$ to isolate the unrepresented residual:
\begin{equation}
\label{eq:s3-d18}
\boldsymbol{R}_i=(\boldsymbol{I}-\widetilde{\boldsymbol{U}}_i\widetilde{\boldsymbol{U}}_i^T)\boldsymbol{W}_i.
\end{equation}

\item \textbf{Enrich the core basis if necessary.} If $\|\boldsymbol{R}_i\|_F>\varepsilon$, augment $\widetilde{\boldsymbol{U}}_i$ with an orthonormal basis for the column space of $\boldsymbol{R}_i$:
\[
\boldsymbol{U}_i^+=[\widetilde{\boldsymbol{U}}_i\;\boldsymbol{P}_i],\qquad \boldsymbol{P}_i=\operatorname{orth}(\boldsymbol{R}_i),\qquad
\boldsymbol{P}_i^T\boldsymbol{P}_i=\boldsymbol{I},\qquad \widetilde{\boldsymbol{U}}_i^T\boldsymbol{P}_i=\boldsymbol{0},
\]
where $\operatorname{orth}(\boldsymbol{R}_i)$ can be implemented via a thin QR factorization. Otherwise, retain $\boldsymbol{U}_i^+=\widetilde{\boldsymbol{U}}_i$. Reshaping $\boldsymbol{U}_i^+$ yields the updated core. This step retains the full residual column space whenever enrichment is triggered, with $\varepsilon$ acting as an absolute local tolerance.

\item \textbf{Propagate the coefficients.} Express $\boldsymbol{W}_i$ in the updated basis through $\boldsymbol{H}_i=(\boldsymbol{U}_i^+)^T\boldsymbol{W}_i$. For $i<d$, contract $\boldsymbol{H}_i$ with the next incoming core to carry the represented information into the next mode:
\begin{equation}
\label{eq:s3-d21}
\mathcal{W}_{i+1}(a,i_{i+1},b)
=\sum_{\gamma=1}^{s_i}\boldsymbol{H}_i(a,\gamma)\mathcal{X}_{i+1}(\gamma,i_{i+1},b).
\end{equation}
Set $\boldsymbol{W}_{i+1}=\operatorname{unfold}_L(\mathcal{W}_{i+1})$ and continue the scan. At $i=d$, $\boldsymbol{H}_d$ yields the coefficient vector of the new snapshot, and we append it as a new column to the coefficient core $\mathcal{U}_{d+1}$.
\end{enumerate}
Optional TT rounding can be  applied to further control the rank growth. Its contribution to the approximation error is included in Section~\ref{sec:4}.

Under the orthogonalization assumptions in Section~\ref{sec:3-1}, Theorem~4.1 of~\cite{shang2026incremental} guarantees that, without TT rounding, the final accumulated tensor approximates each  snapshot with a Frobenius-norm error of at most $\sqrt{d}\varepsilon$. 

Finally, we define the reduced basis
$\boldsymbol{V}\in\mathbb{R}^{(q_1\cdots q_d)\times r_u}$ in TT form by
\[
\boldsymbol{V}((i_1,\ldots,i_d),\alpha)
=
\mathcal{U}_1^{(M)}(1,i_1,:)\,
\mathcal{U}_2^{(M)}(:,i_2,:)\cdots
\mathcal{U}_d^{(M)}(:,i_d,\alpha),
\qquad \alpha=1,\ldots,r_u,
\]
where $r_u$ is the final TT rank connecting the physical cores
to the coefficient core.

\textbf{Projection of affine tensorized operators.}
For the affine linear special case, consider the decomposition~\eqref{eq:s2-d9} where each component has Kronecker structure (with first-index-fastest ordering):
\begin{equation}
\label{eq:s3-affine-tensor-operator}
\boldsymbol{A}_q=\boldsymbol{A}_{q,d}\otimes\cdots\otimes \boldsymbol{A}_{q,1},
\qquad \boldsymbol{A}_{q,k}\in\mathbb R^{n_k\times n_k},
\end{equation}
giving the reduced operator
\begin{equation}
\label{eq:s3-affine-reduced-operator}
\boldsymbol{A}_r(\boldsymbol{\mu})=\sum_{q=1}^{Q_A}\theta_q^A(\boldsymbol{\mu})\boldsymbol{A}_{r,q},
\qquad \boldsymbol{A}_{r,q}=\boldsymbol{V}^T\boldsymbol{A}_q\boldsymbol{V}.
\end{equation}
For finalized cores $\mathcal{U}_k\in\mathbb R^{r_{k-1}\times n_k\times r_k}$ ($r_0=1$, $r_d=r_u$) with left unfoldings $\widetilde{\boldsymbol{U}}_k=\operatorname{unfold}_L(\mathcal{U}_k)$, initialize $\boldsymbol{E}_{q,0}=[1]$ and compute recursively for $k=1,\ldots,d$:
\begin{equation}
\label{eq:s3-affine-core-projection}
\boldsymbol{E}_{q,k}
=\widetilde{\boldsymbol{U}}_k^T
\bigl(\boldsymbol{A}_{q,k}\otimes \boldsymbol{E}_{q,k-1}\bigr)\widetilde{\boldsymbol{U}}_k \in\mathbb R^{r_k\times r_k},
\end{equation}
yielding $\boldsymbol{A}_{r,q}=\boldsymbol{E}_{q,d}$. Applying $\boldsymbol{A}_{q,k}$ to the physical fibers of $\mathcal{U}_k$ and multiplying slices on the left by $\boldsymbol{E}_{q,k-1}$ avoids forming Kronecker products or the full tensor reconstruction of the reduced basis while exploiting sparsity in $\boldsymbol{A}_{q,k}$.

Similarly, for a right-hand side $\boldsymbol{f}_q$ with TT cores $\mathcal{F}_{q,k}\in\mathbb R^{s_{q,k-1}\times n_k\times s_{q,k}}$ ($s_{q,0}=s_{q,d}=1$) and $\boldsymbol{H}_{q,0}=[1]$, define the slices $\boldsymbol{U}_k(i)=\mathcal{U}_k(:,i,:)$ and $\boldsymbol{F}_{q,k}(i)=\mathcal{F}_{q,k}(:,i,:)$, and compute
\begin{equation}
\label{eq:s3-affine-rhs-projection}
\begin{aligned}
\boldsymbol{H}_{q,k}
&=\sum_{i=1}^{n_k}
\boldsymbol{U}_k(i)^T\boldsymbol{H}_{q,k-1}\boldsymbol{F}_{q,k}(i)
\in\mathbb R^{r_k\times s_{q,k}},
\qquad k=1,\ldots,d,\\
\boldsymbol{f}_{r,q}
&=\boldsymbol{H}_{q,d}(:,1)=\boldsymbol{V}^T\boldsymbol{f}_q,
\end{aligned}
\end{equation}
where $s_{q,k}=1$ corresponds to a separable term. All $\boldsymbol{A}_{r,q}$ and $\boldsymbol{f}_{r,q}$ are precomputed offline for online assembly via~\eqref{eq:s3-affine-reduced-operator} and $\boldsymbol{f}_r(\boldsymbol{\mu})=\sum_{q=1}^{Q_f}\theta_q^f(\boldsymbol{\mu})\boldsymbol{f}_{r,q}$.

\subsection{Offline Hyper Reduction through TT-DEIM}
\label{sec:3-3}

To identify the contributions requiring hyper reduction, decompose the residual as
\begin{equation}
\label{eq:s3-residual-decomposition}
\begin{aligned}
\boldsymbol{r}(\boldsymbol{u};\boldsymbol{\mu})
&=\boldsymbol{A}_{\mathrm{aff}}(\boldsymbol{\mu})\boldsymbol{u}-\boldsymbol{f}_{\mathrm{aff}}(\boldsymbol{\mu})
+\boldsymbol{r}_{\mathrm{non}}(\boldsymbol{u};\boldsymbol{\mu}),\\
 \boldsymbol{r}_{\mathrm{non}}(\boldsymbol{u};\boldsymbol{\mu})
&=\sum_{\ell=1}^{L}\boldsymbol{g}_\ell(\boldsymbol{u};\boldsymbol{\mu}),\qquad \boldsymbol{g}_\ell(\boldsymbol{u};\boldsymbol{\mu})\in\mathbb R^N.
\end{aligned}
\end{equation}
Here $\boldsymbol{A}_{\mathrm{aff}}$ and $\boldsymbol{f}_{\mathrm{aff}}$ collect the available affine linear components, whose reduced projections are precomputed as in Section~\ref{sec:3-2}. The remainder $ \boldsymbol{r}_{\mathrm{non}}$ contains a finite sum of nonlinear or non-affine terms. 

Building on nested TT-DEIM~\cite{dektor2025collocation} and the related TT-MDEIM construction~\cite{mueller2026ttrom}, we approximate each contribution separately, with its own collateral basis, interpolation indices, and reduced quantities. In the following termwise construction, we omit the subscript $\ell$ and write $\boldsymbol{g}$ for a generic interpolation target.  The main step  here is the selection of interpolation indices directly from the basis cores, which allows the expansion coefficients to be recovered without evaluating the complete target online.

\textbf{Snapshot generation and basis construction.}
Let $\boldsymbol{g}(\boldsymbol{u};\boldsymbol{\mu})\in\mathbb R^{N_g}$ denote a term to be interpolated, viewed as the vectorization of a tensor with mode sizes $q_1,\ldots,q_d$ and $N_g=\prod_{k=1}^d q_k$. At selected training states and parameters, we collect the snapshots
\[
\boldsymbol{g}^{(j)}=\boldsymbol{g}(\boldsymbol{u}^{(j)};\boldsymbol{\mu}^{(j)}),\qquad j=1,\ldots,J_g.
\]
Each snapshot is expressed in TT format directly if available, or compressed via TT-CUR-DEIM ~\cite{ghahremani2024cross}  (see Appendix~\ref{app:auxiliary-cross}).

After right-orthogonalizing these TT snapshots, we apply the incremental construction of Section~\ref{sec:3-2}. Excluding the sample core, the full-tensor matricization of physical cores define an orthonormal DEIM basis $\boldsymbol{U}_g\in\mathbb R^{N_g\times m}$. We approximate $\boldsymbol{g}$ by
\begin{equation}
\label{eq:s3-term-expansion}
\widehat{\boldsymbol{g}}(\boldsymbol{u};\boldsymbol{\mu})=\boldsymbol{U}_g\boldsymbol{c}_g(\boldsymbol{u};\boldsymbol{\mu}),
\end{equation}
where $\boldsymbol{c}_g(\boldsymbol{u};\boldsymbol{\mu})$ is recovered online by interpolation from selected entries of $\boldsymbol{g}$. 

\textbf{Nested TT-DEIM selection.} We seek the TT representation of
 $\boldsymbol{U}_g\in\mathbb R^{N_g\times m}$, namely $\mathcal{U}_g\in\mathbb R^{q_1\times\cdots\times q_d\times m}$, represented by the physical TT cores
\begin{equation}
\label{eq:s3-d23}
\begin{aligned}
\mathcal{G}_k&\in\mathbb R^{r_{k-1}\times q_k\times r_k},\qquad r_0=1,\quad r_d=m,\\
\mathcal{U}_g(i_1,\ldots,i_d,b)
&=\bigl[\boldsymbol{G}_1(i_1)\cdots \boldsymbol{G}_d(i_d)\bigr]_{1,b},\qquad b=1,\ldots,m,
\end{aligned}
\end{equation}
where $\boldsymbol{G}_k(i_k)=\mathcal{G}_k(:,i_k,:)$ is an $r_{k-1}\times r_k$ matrix. The last core has size $r_{d-1}\times q_d\times m$, with an open right index labeling the $m$ basis vectors.

We seek $m$ physical multi-indices $\mathcal{I}_d$ such that $\boldsymbol{M}_d=\boldsymbol{U}_g[\mathcal{I}_d,:]\in\mathbb R^{m\times m}$ is nonsingular. Nested TT-DEIM~\cite{dektor2025collocation,mueller2026ttrom} selects them through local core contractions without expanding the basis:
\begin{enumerate}
\item \textbf{Initialize.} Set $\mathcal{I}_0=\{()\}$ and $\boldsymbol{M}_0=[1]$. Perform Steps 2--3 for $k=1,\ldots,d$.

\item \textbf{Extend the prefixes and select rows.} The matrix $\boldsymbol{M}_{k-1}\in\mathbb R^{r_{k-1}\times r_{k-1}}$ stores the interface values at the selected prefixes $\mathcal{I}_{k-1}$. Form the local matrix and select $r_k$ rows by DEIM:
\begin{equation}
\label{eq:s3-d24}
\begin{aligned}
\boldsymbol{B}_k&=\operatorname{unfold}_L(\boldsymbol{M}_{k-1}\times_1\mathcal{G}_k)
\in\mathbb R^{(r_{k-1}q_k)\times r_k},\\
\boldsymbol{\ell}_k&=\operatorname{DEIM}(\boldsymbol{B}_k),\qquad \boldsymbol{M}_k=\boldsymbol{B}_k[\boldsymbol{\ell}_k,:].
\end{aligned}
\end{equation}

\item \textbf{Decode the selected indices.} Using first-index-fastest row ordering, write each selected row as $\boldsymbol{\ell}_k(j)=\alpha_j+r_{k-1}(i_{k,j}-1)$ and set
\[
\mathcal{I}_k(j)=\bigl(\mathcal{I}_{k-1}(\alpha_j),i_{k,j}\bigr),\qquad j=1,\ldots,r_k.
\]

\item \textbf{Finalize the interpolation matrix.} Since $r_d=m$, the final index set $\mathcal{I}_d$ contains $m$ entries and
\[
\boldsymbol{M}_d=\boldsymbol{U}_g[\mathcal{I}_d,:]\in\mathbb R^{m\times m}.
\]
\end{enumerate}
Algorithm~\ref{alg:5} summarizes the recursion.

\begin{algorithm}[tbp]
\caption{Nested TT-DEIM sampling of a basis with an open terminal index}
\label{alg:5}
\footnotesize
\begin{algorithmic}[1]
\Require Basis cores $\mathcal{G}_k\in\mathbb R^{r_{k-1}\times q_k\times r_k}$ of $\boldsymbol{U}_g$ as in~\eqref{eq:s3-d23}, with full-column-rank left unfoldings; $r_0=1$, $r_d=m$
\Ensure $\mathcal{I}_d,\{\mathcal{I}_k\}_{k=1}^{d},\boldsymbol{M}_d=\boldsymbol{U}_g[\mathcal{I}_d,:]$
\State $\mathcal{I}_0\gets\{()\}$; $\boldsymbol{M}_0\gets[1]$
\For{$k=1,\ldots,d$}
  \State $\boldsymbol{B}_k\gets\operatorname{unfold}_L(\boldsymbol{M}_{k-1}\times_1\mathcal{G}_k)$ \Comment{$(r_{k-1}q_k)\times r_k$}
  \State $\boldsymbol{\ell}_k\gets\algcall{DEIM}(\boldsymbol{B}_k)$; $\boldsymbol{M}_k\gets \boldsymbol{B}_k[\boldsymbol{\ell}_k,:]$
  \For{$j=1,\ldots,r_k$}
    \State $(\alpha_j,i_{k,j})\gets\algcall{DecodeIndex}(\boldsymbol{\ell}_k(j);r_{k-1},q_k)$
    \State $\mathcal{I}_k(j)\gets(\mathcal{I}_{k-1}(\alpha_j),i_{k,j})$
  \EndFor
\EndFor
\State \Return $\mathcal{I}_d,\{\mathcal{I}_k\}_{k=1}^{d},\boldsymbol{M}_d$
\end{algorithmic}
\end{algorithm}

Nestedness ensures that the recursively computed interpolation matrix agrees exactly with the sampled full basis. Since $\mathcal{I}_k\subseteq\mathcal{I}_{k-1}\times[q_k]$, the interface matrices satisfy
\[
\boldsymbol{M}_k(j,:)=\boldsymbol{G}_1(i_1)\cdots \boldsymbol{G}_k(i_k),\qquad
(i_1,\ldots,i_k)=\mathcal{I}_k(j).
\]
At the final mode, $\boldsymbol{M}_d=\boldsymbol{U}_g[\mathcal{I}_d,:]$. Therefore, if $\boldsymbol{M}_d$ is nonsingular, the coefficients in~\eqref{eq:s3-term-expansion} yield exact reconstruction at the interpolation points:
\[
\widehat{\boldsymbol{g}}(\boldsymbol{u};\boldsymbol{\mu})[\mathcal{I}_d]
=\boldsymbol{U}_g[\mathcal{I}_d,:]\boldsymbol{c}_g(\boldsymbol{u};\boldsymbol{\mu})
=\boldsymbol{M}_d\boldsymbol{c}_g(\boldsymbol{u};\boldsymbol{\mu})
=\boldsymbol{g}(\boldsymbol{u};\boldsymbol{\mu})[\mathcal{I}_d].
\]

Since the core left unfoldings have full column rank, nonsingular $\boldsymbol{M}_{k-1}$ implies full-column-rank $\boldsymbol{B}_k$, and DEIM selects nonsingular $\boldsymbol{M}_k$. Starting from $\boldsymbol{M}_0=[1]$, induction gives nonsingular $\boldsymbol{M}_d$~\cite{dektor2025collocation}; its conditioning is discussed in Section~\ref{sec:4}.

\textbf{Offline precomputation and interpolation.}
Given the TT cores $\{\mathcal{G}_k\}_{k=1}^d$ for the targeted non-affine or nonlinear term and the indices $\mathcal{I}_d$ from Algorithm~\ref{alg:5}, prepare and apply the interpolant through core contractions and selected-entry evaluations:
\begin{enumerate}
\item \textbf{Factor the sampled interface (offline).} Store a factorization of the $m\times m$ nested-selection matrix $\boldsymbol{M}_d$, along with the basis cores and fixed indices $\mathcal{I}_d$. Throughout, $\boldsymbol{P}_g^T\boldsymbol{x}$ denotes extracting the entries at $\mathcal{I}_d$ without assembling a selection matrix.

\item \textbf{Precompute the reduced maps (offline).} For a vector term $\boldsymbol{g}$ on the state grid, equation~\eqref{eq:s2-deim-preprojection} reads, in the current notation,
\[
\begin{aligned}
\boldsymbol{V}^T\widehat{\boldsymbol{g}}=(\boldsymbol{V}^T\boldsymbol{U}_g)(\boldsymbol{P}_g^T\boldsymbol{U}_g)^{-1}\boldsymbol{P}_g^T\boldsymbol{g}=(\boldsymbol{V}^T\boldsymbol{U}_g)[\boldsymbol{U}_g(\mathcal{I}_d,:)]^{-1}\boldsymbol{g}(\mathcal{I}_d).
\end{aligned}
\]
We therefore precompute
\[
\boldsymbol{B}_{r,g}=\boldsymbol{V}^T\boldsymbol{U}_g\in\mathbb R^{r_u\times m}\quad\text{and}\quad
\boldsymbol{Z}_g=\boldsymbol{P}_g^T\boldsymbol{V}\in\mathbb R^{m\times r_u}.
\]
Compute $\boldsymbol{B}_{r,g}$ via the core contraction~\eqref{eq:s3-affine-rhs-projection}, retaining the $m$-dimensional open terminal index of $\boldsymbol{U}_g$. Obtain $\boldsymbol{Z}_g$ by evaluating the solution-basis core slices $\boldsymbol{U}_k(i)=\mathcal{U}_k(:,i,:)$ at the selected multi-indices:
\[
\boldsymbol{Z}_g(j,:)=\boldsymbol{U}_1(i_{1,j})\cdots \boldsymbol{U}_d(i_{d,j}),\qquad
(i_{1,j},\ldots,i_{d,j})=\mathcal{I}_d(j).
\]

\end{enumerate}

\subsection{Online Solve}
\label{sec:3-4}

At a new parameter $\boldsymbol{\mu}$, using the fixed offline quantities from Section~\ref{sec:3-3} and the decomposition~\eqref{eq:s3-residual-decomposition}, we solve the hyper-reduced system for $\boldsymbol{c}_u(\boldsymbol{\mu})$:
\[
\widehat{\boldsymbol{R}}_r(\boldsymbol{c}_u;\boldsymbol{\mu}) = \boldsymbol{A}_{\mathrm{aff},r}(\boldsymbol{\mu})\boldsymbol{c}_u - \boldsymbol{f}_{\mathrm{aff},r}(\boldsymbol{\mu}) + \sum_{\ell=1}^{L} \boldsymbol{V}^T\widehat{\boldsymbol{g}}_\ell(\boldsymbol{V}\boldsymbol{c}_u;\boldsymbol{\mu}) = \boldsymbol{0},
\]
where $\boldsymbol{A}_{\mathrm{aff},r}=\boldsymbol{V}^T\boldsymbol{A}_{\mathrm{aff}}\boldsymbol{V}$ and $\boldsymbol{f}_{\mathrm{aff},r}=\boldsymbol{V}^T\boldsymbol{f}_{\mathrm{aff}}$ are assembled from precomputed affine projections, and each $\widehat{\boldsymbol{g}}_\ell$ is constructed via contractions and inerpolations. Starting from an initial iterate $\boldsymbol{c}_u^{(0)}$, the online solve proceeds as follows:
\begin{enumerate}
\item \textbf{Sample TT cores, interpolate, and assemble the reduced residual.}
For each $\ell=1,\ldots,L$, we approximate the targeted  non-affine or nonlinear term at the selected multi-indices $\mathcal{I}_{d,\ell}$ using the solution-basis cores:
$$\begin{gathered}
(\boldsymbol{V}\boldsymbol{c}_u)(i_{1,j},\ldots,i_{d,j}) = \boldsymbol{U}_1(i_{1,j})\cdots \boldsymbol{U}_d(i_{d,j})\boldsymbol{c}_u = (\boldsymbol{Z}_{g_\ell}\boldsymbol{c}_u)_j,\\
(i_{1,j},\ldots,i_{d,j})=\mathcal{I}_{d,\ell}(j).
\end{gathered}$$
Using these sampled  values, along with any additional core contractions required by the sampled term, evaluate the targeted term as
$$\boldsymbol{y}_{g_\ell} = \boldsymbol{g}_\ell(\boldsymbol{V}\boldsymbol{c}_u;\boldsymbol{\mu})[\mathcal{I}_{d,\ell}]$$
without assembling the full tensor. Next, recover the interpolation coefficients $\boldsymbol{c}_{g_\ell}$ by solving the small linear system

$$\boldsymbol{M}_{d,\ell}\boldsymbol{c}_{g_\ell}=\boldsymbol{y}_{g_\ell}$$

via the stored offline factorization, and assemble the reduced residual:

$$\widehat{\boldsymbol{R}}_r(\boldsymbol{c}_u;\boldsymbol{\mu}) = \boldsymbol{A}_{\mathrm{aff},r}(\boldsymbol{\mu})\boldsymbol{c}_u-\boldsymbol{f}_{\mathrm{aff},r}(\boldsymbol{\mu}) + \sum_{\ell=1}^{L}\boldsymbol{B}_{r,g_\ell}\boldsymbol{c}_{g_\ell}.$$
Consequently, this online stage relies only on pointwise evaluations, TT-core contractions, and small-scale operations.

\item \textbf{Solve for the reduced coefficients.} Determine $\boldsymbol{c}_u\in\mathbb R^{r_u}$ by the nonlinear reduced system
\[
\boldsymbol{A}_{\mathrm{aff},r}(\boldsymbol{\mu})\boldsymbol{c}_u
+\sum_{\ell=1}^{L}\boldsymbol{B}_{r,g_\ell}\boldsymbol{c}_{g_\ell}(\boldsymbol{V}\boldsymbol{c}_u;\boldsymbol{\mu})
=\boldsymbol{f}_{\mathrm{aff},r}(\boldsymbol{\mu}).
\]
This system can be solved using Newton's method or nother suitable nonlinear solvers, with the sampled terms and interpolation coefficients reevaluated through Step~1 if necessary.

\item \textbf{Recover the TT solution.} Once converged, form $\boldsymbol{u}_H = \boldsymbol{V}\boldsymbol{c}_u$ in TT format by contracting $\boldsymbol{c}_u$ with the TT formatted reduced basis $\boldsymbol{V}$, without the need of full tensor expansion.
\end{enumerate}

\subsection{Offline Storage and Construction Cost}
\label{sec:3-5}

We summarize the costs of the constituent algorithms within the offline procedure.

\textbf{Storage.} For a data stream of \(J\) snapshots with \(N = \prod_{k=1}^d q_k\) and \(q = \max_k q_k\), let \(R\) bound all TT ranks throughout construction, including the final basis dimension \(m \le R\). Compared with standard POD--DEIM, which stores the snapshots and basis in vectorized form, the offline storage requirements are
\[
\begin{aligned}
\text{TT storage} &= \mathcal{O}\!\left(dqR^2 + JR + m^2 + dm\right),\\
\text{Full-snapshot storage} &= \mathcal{O}\!\left(N_g(J+m) + m^2 + dm\right),
\end{aligned}
\]
where \(dqR^2\) and \(JR\) account for the physical TT cores and sample coefficients~\cite[Section~4.2]{shang2026incremental}, while \(m^2\) and \(dm\) store the interpolation factor and sampling multi-indices.

\textbf{Computational cost.} Under the orthogonality assumptions of Section~\ref{sec:3-1}, incremental reduced basis construction performs residual projection, thin QR factorization for basis enrichment, and coefficient propagation on core unfoldings of size at most $qR\times R$. Each operation costs $\mathcal{O}(qR^3)$ per mode and snapshot~\cite[Section~4.2]{shang2026incremental}. Nested TT-DEIM applies DEIM to a local matrix of size at most $qR\times R$ at each mode~\cite[Section~5.1]{dektor2025collocation}. With fresh factorizations in greedy DEIM and the interpolation-matrix factorization of the sampling core included, the total costs are
\begin{equation}
\label{eq:s3-offline-compression}
\begin{aligned}
&\text{Incremental TT basis construction} = \mathcal{O}\!\left(JdqR^3\right),\\
&\text{Nested TT-DEIM} = \mathcal{O}\!\left(d(qR^3+R^4)+m^3\right).
\end{aligned}
\end{equation}
These compression estimates exclude snapshot generation, optional TT rounding, additional orthogonalization sweeps, and offline reduced-map precomputation.

Dense full tensor POD costs $\mathcal{O}(N_gJ\min\{N_g,J\})$, while DEIM on an explicit basis costs $\mathcal{O}(N_gm^2+m^4)$ under the same fresh-factorization convention. For equal mode sizes $q_k=q>1$, the full-order dimension is $N_g=q^d$. When $J\le N_g$, the dense costs therefore become
\[
\begin{aligned}
&\text{Dense SVD/POD} = O(q^dJ^2),\\
&\text{Explicit-basis DEIM} = O(q^dm^2+m^4).
\end{aligned}
\]
If degrees of freedom in each dimension is roughly the same, for fixed $q$, $J$, and $m$, these costs grow exponentially with $d$. In contrast, the TT arithmetic costs in~\eqref{eq:s3-offline-compression} grow linearly with $d$, if the TT rank bound $R$ remains independent of $d$.

\section{Error Bounds for Incremental TT-DEIM}
\label{sec:4}

This section establishes a uniform training interpolation bound for a nonlinear or non-affine term \(\boldsymbol{g}(\boldsymbol{u};\boldsymbol{\mu})\in\mathbb R^{N_g}\) approximated via LR-iTT and nested TT-DEIM (Corollary~\ref{cor:s4-training-tt-deim}). The bound separates the approximation error of incrementally building the TT reduced basis from TT snapshots and the interpolation error due to TT-DEIM.

Our proof adapts~\cite[Theorem~3]{mueller2026ttrom}, replacing its TT-SVD estimate with an incremental LR-iTT bound. Specifically, Lemma~\ref{lem:s4-lritt} bounds the basis projection error on the incoming TT snapshots via~\cite[Theorem~4.1]{shang2026incremental}, while Lemma~\ref{lem:s4-tt-deim} applies~\cite[Theorem~5.1]{dektor2025collocation} to guarantee nonsingularity of the sampled matrices in Algorithm~\ref{alg:5} and obtain the DEIM interpolation estimate. Combining both bounds via the triangle inequality yields Corollary~\ref{cor:s4-training-tt-deim}.

\begin{lemma}[LR-iTT projection error ]\cite[Theorem~4.1]{shang2026incremental}
\label{lem:s4-lritt}
For each training pair \((\boldsymbol{u}^{(j)},\boldsymbol{\mu}^{(j)})\), \(j=1,\ldots,J_g\), let \(\boldsymbol{g}^{(j)}=\boldsymbol{g}(\boldsymbol{u}^{(j)};\boldsymbol{\mu}^{(j)})\) be the vectorized full-tensor snapshot. Let \(\mathcal{G}_1^{(j)},\ldots,\mathcal{G}_d^{(j)}\) be the TT cores of its approximation supplied to the streaming algorithm, with physical mode sizes \(q_1,\ldots,q_d\) and unit boundary ranks. Define
\[
\boldsymbol{g}_{\mathrm{TT}}^{(j)}
:=\operatorname{vec}\!\left(
\left[\mathcal{G}_1^{(j)}(:,i_1,:)\cdots \mathcal{G}_d^{(j)}(:,i_d,:)\right]_{i_1,\ldots,i_d}
\right)\in\mathbb R^{N_g},
\]
where \(1\le i_k\le q_k\), \(N_g=\prod_{k=1}^d q_k\), and vectorization uses the same index ordering as \(\boldsymbol{g}^{(j)}\). Thus \(\boldsymbol{g}_{\mathrm{TT}}^{(j)}\) is the vectorized incoming TT approximation before incremental basis construction. The algorithm processes its cores directly, without forming this vector.

Construct the orthonormal basis \(\boldsymbol{U}_g\in\mathbb R^{N_g\times m}\) from these TT snapshots by the LR-iTT-F procedure of Section~\ref{sec:3-2}, with absolute local residual tolerance \(\varepsilon_{\mathrm{iTT}}\) at each TT core and the orthogonalization conditions of Section~\ref{sec:3-1}. Assume no optional TT rounding is performed during the incremental updates or after the final update. Then
\begin{equation}
\label{eq:s4-lritt-projection}
\max_{1\le j\le J_g}\|(\boldsymbol{I}-\boldsymbol{U}_g\boldsymbol{U}_g^T)\boldsymbol{g}_{\mathrm{TT}}^{(j)}\|_2
\le\sqrt d\,\varepsilon_{\mathrm{iTT}}.
\end{equation}
\end{lemma}

\begin{lemma}[Nested TT-DEIM interpolation]
\label{lem:s4-tt-deim}
Let \(\boldsymbol{U}_g\in\mathbb R^{N_g\times m}\) have TT cores \(\mathcal{G}_k\in\mathbb R^{r_{k-1}\times q_k\times r_k}\) as in~\eqref{eq:s3-d23}, with \(r_0=1\) and \(r_d=m\). Assume the physical cores are left-orthogonal, as produced by LR-iTT:
\[
\operatorname{unfold}_L(\mathcal{G}_k)^T\operatorname{unfold}_L(\mathcal{G}_k)=\boldsymbol{I}_{r_k},
\qquad k=1,\ldots,d.
\]
Let \(\mathcal{I}_k\) and \(\boldsymbol{M}_k\in\mathbb R^{r_k\times r_k}\) be the index sets and sampled matrices computed by Algorithm~\ref{alg:5}, and let \(\boldsymbol{P}_g^T\boldsymbol{x}=\boldsymbol{x}[\mathcal{I}_d]\). Then each \(\boldsymbol{M}_k\) is nonsingular and \(\boldsymbol{M}_d=\boldsymbol{P}_g^T\boldsymbol{U}_g\). For every \(\boldsymbol{g}\in\mathbb R^{N_g}\), the TT-DEIM interpolant satisfies
\begin{equation}
\label{eq:s4-d15}
\begin{aligned}
&\|(\boldsymbol{I}-\boldsymbol{U}_g(\boldsymbol{P}_g^T\boldsymbol{U}_g)^{-1}\boldsymbol{P}_g^T)\boldsymbol{g}\|_2\le\|(\boldsymbol{P}_g^T\boldsymbol{U}_g)^{-1}\|_2\,\|(\boldsymbol{I}-\boldsymbol{U}_g\boldsymbol{U}_g^T)\boldsymbol{g}\|_2.
\end{aligned}
\end{equation}
\end{lemma}

\begin{proof}
The sampled matrices \(\boldsymbol{M}_k\) contain contractions of the first \(k\) basis cores at the selected prefixes:
\[
\boldsymbol{M}_k(j,:)=\boldsymbol{G}_1(i_1)\cdots \boldsymbol{G}_k(i_k),
\qquad (i_1,\ldots,i_k)=\mathcal{I}_k(j).
\]
These are the sampled matrices in Dektor~\cite[Eqs.~(47) and~(59)]{dektor2025collocation}, with the basis index viewed as an additional terminal mode. The assumed left orthogonality gives \(\operatorname{rank}(\operatorname{unfold}_L(\mathcal{G}_k))=r_k\) for every \(k\), which is the full-column-rank condition used in the left-sweep argument of Dektor~\cite[Theorem~5.1]{dektor2025collocation}. Its nonsingularity conclusion therefore applies to every \(\boldsymbol{M}_k\) computed by Algorithm~\ref{alg:5}. In particular,
\[
\boldsymbol{M}_d=\boldsymbol{U}_g[\mathcal{I}_d,:]=\boldsymbol{P}_g^T\boldsymbol{U}_g.
\]
Left orthogonality of the physical cores also gives \(\boldsymbol{U}_g^T\boldsymbol{U}_g=\boldsymbol{I}_m\).
For interpolation, \((\boldsymbol{I}-\boldsymbol{U}_g(\boldsymbol{P}_g^T\boldsymbol{U}_g)^{-1}\boldsymbol{P}_g^T)\boldsymbol{U}_g=\boldsymbol{0}\). Hence
\[
\begin{aligned}
\|(\boldsymbol{I}-\boldsymbol{U}_g(\boldsymbol{P}_g^T\boldsymbol{U}_g)^{-1}\boldsymbol{P}_g^T)\boldsymbol{g}\|_2
&=\|(\boldsymbol{I}-\boldsymbol{U}_g(\boldsymbol{P}_g^T\boldsymbol{U}_g)^{-1}\boldsymbol{P}_g^T)(\boldsymbol{I}-\boldsymbol{U}_g\boldsymbol{U}_g^T)\boldsymbol{g}\|_2\\
&\le\|(\boldsymbol{P}_g^T\boldsymbol{U}_g)^{-1}\|_2\,\|(\boldsymbol{I}-\boldsymbol{U}_g\boldsymbol{U}_g^T)\boldsymbol{g}\|_2,
\end{aligned}
\]
where the last inequality follows from orthonormality of \(\boldsymbol{U}_g\) and the standard DEIM projector bound~\cite{chaturantabut2010deim}. This proves~\eqref{eq:s4-d15}.
\end{proof}

\begin{theorem}[Training TT-DEIM error]
\label{cor:s4-training-tt-deim}
Under the assumptions of Lemmas~\ref{lem:s4-lritt} and~\ref{lem:s4-tt-deim}, suppose the supplied TT snapshots satisfy the uniform bound
\[
\max_{1\le j\le J_g}\|\boldsymbol{g}^{(j)}-\boldsymbol{g}_{\mathrm{TT}}^{(j)}\|_2
\le\varepsilon_{\mathrm{snap}}.
\]
Then the training interpolation error satisfies
\begin{equation}
\label{eq:s4-training-tt-deim}
\begin{aligned}
&\max_{1\le j\le J_g}\|(\boldsymbol{I}-\boldsymbol{U}_g(\boldsymbol{P}_g^T\boldsymbol{U}_g)^{-1}\boldsymbol{P}_g^T)\boldsymbol{g}^{(j)}\|_2\le\|(\boldsymbol{P}_g^T\boldsymbol{U}_g)^{-1}\|_2\bigl(\varepsilon_{\mathrm{snap}}+\sqrt d\,\varepsilon_{\mathrm{iTT}}\bigr).
\end{aligned}
\end{equation}
\end{theorem}

\begin{proof}
Contractivity of orthogonal projection and Lemma~\ref{lem:s4-lritt} give, for every \(j\),
\begin{equation}
\label{eq:s4-d14}
\begin{aligned}
\|(\boldsymbol{I}-\boldsymbol{U}_g\boldsymbol{U}_g^T)\boldsymbol{g}^{(j)}\|_2
&\le\|\boldsymbol{g}^{(j)}-\boldsymbol{g}_{\mathrm{TT}}^{(j)}\|_2+\|(\boldsymbol{I}-\boldsymbol{U}_g\boldsymbol{U}_g^T)\boldsymbol{g}_{\mathrm{TT}}^{(j)}\|_2\\
&\le\varepsilon_{\mathrm{snap}}+\sqrt d\,\varepsilon_{\mathrm{iTT}}.
\end{aligned}
\end{equation}
Applying Lemma~\ref{lem:s4-tt-deim} and taking the maximum over \(j\) proves~\eqref{eq:s4-training-tt-deim}.
\end{proof}

\section{Numerical Experiments}
\label{sec:5}

We evaluate the  hyper-reduced TT-ROM built from streaming TT data on four problems that introduce non-affine coefficients, nonlinearity, time dependence, and space--angle coupling. Section~\ref{sec:5-1} describes the common data sources and evaluation protocol. Sections~\ref{sec:5-2}--\ref{sec:5-5} present the diffusion, Bratu, Allen--Cahn, and radiative transfer examples, respectively. Section~\ref{sec:5-6} summarizes the accuracy, solution storage, and online efficiency across the four examples and assesses when the online savings offset the additional offline cost of hyper reduction.

\FloatBarrier
\subsection{Data Sources and Experimental Setup}
\label{sec:5-1}

We compare the TT-ROM with TT-DEIM-based hyper reduction against the TT Galerkin ROM without hyper reduction. Both ROMs are built from the same low-rank training data. The corresponding low-rank solver provides an additional baseline for assessing the online computational savings obtained through model reduction.

\textbf{Low-rank snapshot generation.} Diffusion, Bratu, and Allen--Cahn snapshots are generated using the low-rank Anderson acceleration (lrAA) \cite{appelo2026lraa}. It applies Anderson acceleration to matrix fixed-point equations in low-rank format, controls intermediate ranks through adaptive truncation, and uses Cross-DEIM to approximate the fixed-point map from selected matrix entries through an implementation following Appendix \ref{app:auxiliary-cross}. RTE snapshots are generated using the rank-adaptive, sweep-based source iteration with diffusion synthetic acceleration (SI--DSA) ~\cite{guo2026sweep}. Both solvers compute low-rank matrix approximations. In our implementations, both solvers directly approximate the solution in SVD form.

\textbf{Training and test data.} Independent randomized Sobol sequences~\cite{sobol1967distribution} generate the training (five seeds) and test (three seeds) sets, with test data excluded from basis construction and every trained model evaluated on all test points. For diffusion, Bratu, and Allen--Cahn, training seeds contain 64, 64, and 32 points respectively (excluding the center), while each test seed provides 64 points (192 test parameters; 960 test runs per model). For RTE, each training seed includes 32 Sobol points plus the center, and each test seed uses the first five points of an eight-point Sobol sequence (15 test parameters; 75 test runs).

\textbf{Thresholds settings.}  The local rank threshold $\varepsilon_{\mathrm{rank}}$ controls rank increases through the test on the first omitted singular value, while $\varepsilon_{\mathrm{iTT}}$ controls incremental basis updates. Both thresholds are specified for each case.

\textbf{Error and timing protocol.} Steady-state accuracy is measured against the full-order reference solution using the relative Euclidean error:
\begin{equation}
\label{eq:s5-d1}
e_M(\boldsymbol{\mu})=\frac{\|\boldsymbol{u}_M(\boldsymbol{\mu})-\boldsymbol{u}_{\mathrm{FOM}}(\boldsymbol{\mu})\|_2}{\|\boldsymbol{u}_{\mathrm{FOM}}(\boldsymbol{\mu})\|_2},\qquad M\in\{\mathrm{LR},\mathrm G,\mathrm H\}.
\end{equation}
Here LR, G, and H denote the low-rank solver, TT Galerkin ROM, and TT hyper-ROM, respectively; time-dependent errors are defined in Section~\ref{sec:5-4}. Total online runtime for each ROM sums the median times (over five runs after one warm-up) across four stages: input preparation, reduced assembly, reduced solve, and TT reconstruction. All timings use single-threaded Python 3.12.4 on an AMD Ryzen 7 5700U with 16 GB RAM. Across $N_{\mathrm{pair}}$ paired queries with runtimes $T_{\mathrm{G},j}$ and $T_{\mathrm{H},j}$, we report the mean speedup and mean time saving:
\begin{equation}
\label{eq:s5-d2}
S_{\mathrm{HG}}=\frac{1}{N_{\mathrm{pair}}}\sum_{j=1}^{N_{\mathrm{pair}}}\frac{T_{\mathrm G,j}}{T_{\mathrm H,j}},\qquad
\overline{\Delta T}_{\mathrm{on}}=\frac{1}{N_{\mathrm{pair}}}\sum_{j=1}^{N_{\mathrm{pair}}}(T_{\mathrm G,j}-T_{\mathrm H,j}).
\end{equation}

\FloatBarrier
\subsection{non-affine Parametric Diffusion Problem}
\label{sec:5-2}

Consider $-\nabla\cdot(a\nabla u)=f$ on $\Omega=(0,1)^2$, with the boundary condition $u|_{\partial\Omega}=0$, where
\begin{align*}
a(x,y;\boldsymbol{\mu})&=\exp(\mu_1\sin(2\pi x)\sin(2\pi y)+\mu_2\cos(\pi x)\cos(2\pi y)),\\
f(x,y;\boldsymbol{\mu})&=\exp(-(\;(x-\mu_3)^2+(y-\mu_4)^2)/(2(0.12)^2)\;).
\end{align*}
The parameter domain is $[-0.5,0.5]^2\times[0.25,0.75]^2$. Both the diffusion coefficient and the right-hand side have non-affine dependence on the parameters.

The domain is discretized with a second order central difference scheme using \(63\times63\) interior grid points, evaluating diffusion coefficients at grid faces. Sparse LU gives the full-order reference solution. Training uses lrAA with a relative residual tolerance of \(10^{-4}\). The operator-coefficient and right-hand side rank thresholds are both \(10^{-5}\). The TT reduced basis are built using an incremental tolerance of \(10^{-5}\). Across five sets of training data, the solution space dimension is consistently 64, the operator DEIM space dimension is 44-48, and the right-hand-side DEIM space dimension is 64.

\textbf{Accuracy and online efficiency.} On the $63\times63$ grid, the hyper-ROM has a mean relative error of $2.280\times10^{-4}$, compared with $1.676\times10^{-4}$ for the TT Galerkin ROM, and achieves a mean paired speedup of 8.42. 

\begin{figure}[htbp]
\centering
\includegraphics[width=\textwidth]{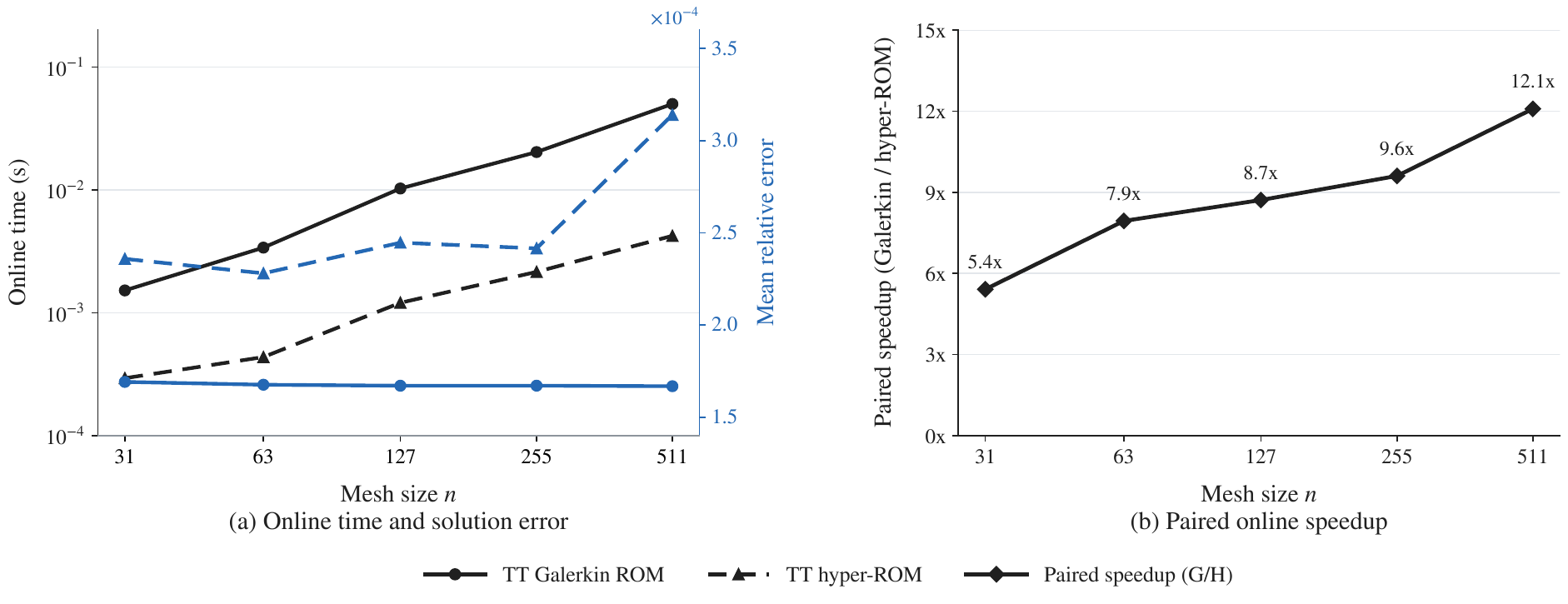}
\caption{Online time, mean relative solution error, and paired speedup for the parametric diffusion problem under mesh refinement. The horizontal axis represents the number of grid points n in each spatial direction, using a logarithmic scale with base 2. The left vertical axis (black line) represents online time, using a logarithmic scale; the right vertical axis (blue line) represents the mean relative solution error, in units of \(10^{-4}\). Solid lines with circular markers represent the TT Galerkin ROM, and dashed lines with triangular markers represent the TT hyper-ROM.}
\label{fig:1}
\end{figure}

\textbf{Effect of grid refinement.} The relative benefit of hyper reduction increases with mesh refinement (Fig.~\ref{fig:1}). As the number of interior grid points per direction increases from 31 to 511, the mean speedup of hyper-reduced TT-ROM over the TT-Galerkin-ROM without hyper reduction  increases from approximately 5.4 to 12.1, while the mean relative error of the hyper-reduced model remains between \(2.3\times10^{-4}\) and \(3.1\times10^{-4}\). This is because the Galerkin model requires an increasing number of coefficient terms to be evaluated and projected as the mesh is refined, whereas the online cost of hyper reduction depends on the number of sampling points and the reduced dimension.

\FloatBarrier
\subsection{Parametric Bratu Problem}
\label{sec:5-3}

We consider a parametric variant of the Bratu problem~\cite{bratu1914equations}: solve $\Delta u+\kappa(x,y;\boldsymbol{\mu})\exp(u)=0$ on $\Omega=(0,1)^2$, subject to homogeneous Dirichlet boundary conditions, where
\[
\kappa=\lambda \exp(\;\mu_1\sin(2\pi x)\sin(\pi y)+\mu_2\cos(\pi x)\sin(2\pi y)\;),
\]
$\lambda\in[0.5,5]$, and $\mu_1,\mu_2\in[-0.4,0.4]$. This problem combines non-affine parameter dependence in $\kappa$ with an exponential nonlinearity in the state variable.

The spatial domain consists of \(63\times63\) interior grid points. The full-order reference solution is obtained using Newton method with sparse LU to solve the Jacobi system in each iteration. The stopping threshold for the lrAA method used in training is \(10^{-6}\). The incremental construction tolerance for both the solution basis and the DEIM basis for the nonlinear term is \(10^{-3}\), and the local rank threshold for nonlinear terms is \(3\times10^{-5}\). The solution space dimensions obtained from the five sets of training data range from 9 to 11, while the dimensions of the DEIM space for the nonlinear term range from 17 to 19; the latter corresponds to the number of online sampling points.

The fixed linear Laplace term is pre-projected offline. The TT-DEIM basis and corresponding interpolation points for the nonlinear term is constructed from \(\kappa\exp(u)\) evaluated at converged training solutions. The resulting nonlinear reduced system is solved  by damped Newton method. Residuals are normalized by the initial reduced residual norm, with a stopping tolerance \(10^{-10}\).

\textbf{Accuracy and online efficiency.} On the $63\times63$ grid, the TT-Galerkin-DEIM achieves a mean paired speedup of 10.39 over the TT-ROM without hyper reduction. Their mean relative errors are $2.825\times10^{-4}$ and $1.905\times10^{-4}$, respectively. The small absolute saving per solve makes the number of repeated queries important to the total-cost comparison in Section~\ref{sec:5-6}.

\begin{figure}[htbp]
\centering
\includegraphics[width=\textwidth]{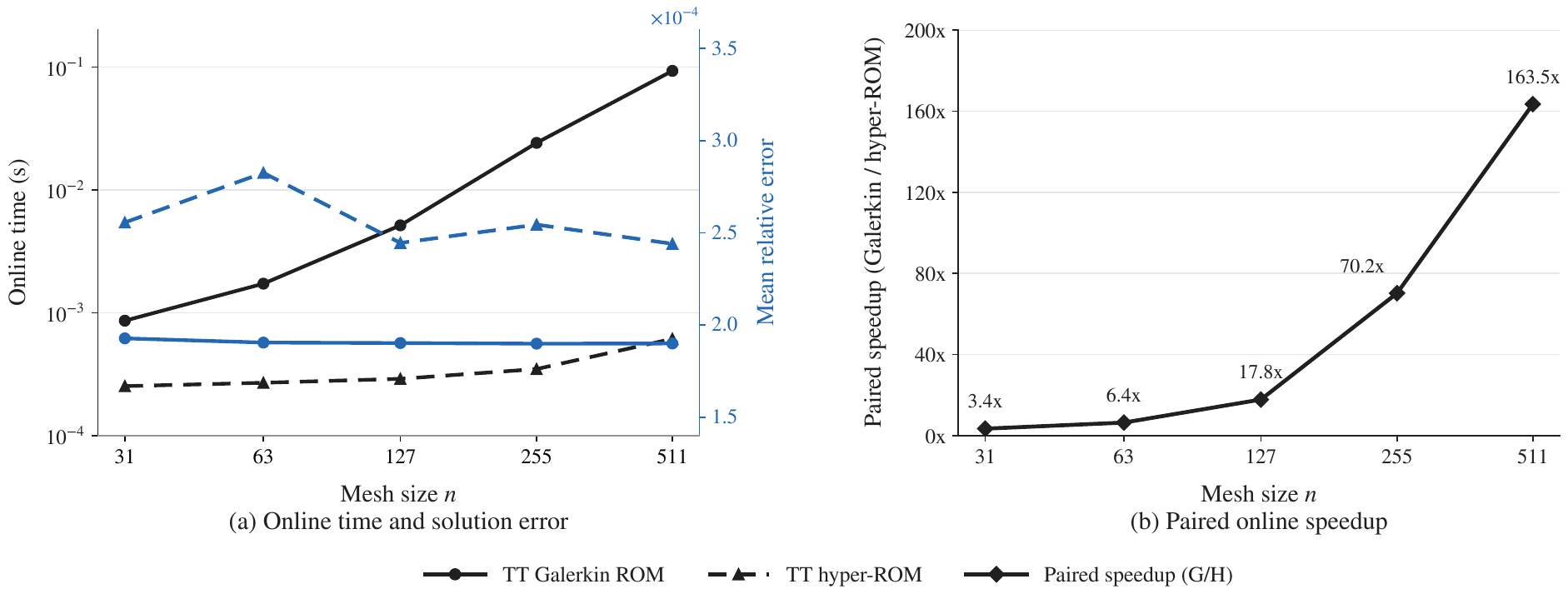}
\caption{Online time, mean relative solution error, and paired speedup for the parametric Bratu problem under grid refinement. The axes, line styles, and markers are defined as in Figure~\ref{fig:1}.}
\label{fig:2}
\end{figure}

\textbf{Effect of grid refinement.} Figure~\ref{fig:2} shows that as the mesh becomes finer, the online time of the Galerkin model increases significantly faster than that of the hyper-reduced model. The speedup gained by hyper reduction increases from about 3.4 to 163.5, while the mean relative errors remain at a similar order of magnitude.

\FloatBarrier
\subsection{Allen--Cahn Equation}
\label{sec:5-4}

To test the framework's applicability to nonlinear time-dependent problems, consider a parametric variant of the Allen--Cahn equation~\cite{allen1979microscopic}

\begin{equation}
\label{eq:s5-d4}
\begin{gathered}
\frac{\partial u}{\partial t}-\nu\Delta u-\kappa(x,y;\mu)(u-u^3)=0,\\
\kappa(x,y;\mu)=e^{(\mu\sin x\sin y)}.
\end{gathered}
\end{equation}
Let $\nu=0.01$, $\mu\in[-0.3,0.3]$, and the spatial domain be $[0,2\pi)^2$; periodic boundary conditions are applied, and parameter dependence is introduced solely through the reaction coefficient $\kappa$. The initial condition is independent of parameters and is defined as
\[
u_0(x,y)=\frac{[e^{(-\tan^2x)}+e^{(-\tan^2y)}]\sin x\sin y}{1+e^{|\csc(\frac x{2})|}+e^{|\csc(\frac y{2})|}},
\]
with continuous extensions at singular points.

A uniform \(256\times256\) grid and backward Euler with \(\Delta t=0.1\) give 50 steps over \([0, 5]\). States are saved every five steps, yielding 11 saved time levels including the initial condition. The reference solver is Newton–GMRES and training uses lrAA with a tolerance of \(10^{-5}\). The solution basis and the DEIM basis for the nonlinear term are constructed from the low-rank states of each training parameter at the saved time levels and their corresponding reaction terms \(\kappa(u-u^3)\) respectively. Both incremental tolerances are \(10^{-4}\) and the local rank threshold is \(5\times10^{-6}\). The resulting solution space has dimensions ranging from 39 to 42, and the DEIM space for the nonlinear term has dimensions ranging from 67 to 71. Online integration uses damped Newton with a fixed DEIM space and fixed sampling indices at every time step.

\textbf{Errors and timing.} For each test parameter, the error is calculated as the equally weighted average of the relative state errors across all saved time steps:

\begin{equation}
\label{eq:s5-d6}
\overline e_M(\mu)=\frac1{11}\sum_{\ell=0}^{10}\frac{\|\boldsymbol{u}_M(t_\ell;\mu)-\boldsymbol{u}_{\mathrm{FOM}}(t_\ell;\mu)\|_2}{\|\boldsymbol{u}_{\mathrm{FOM}}(t_\ell;\mu)\|_2},\qquad t_\ell=0.5\ell.
\end{equation}
The maximum error is calculated across all cases and saved time levels. The online time includes all 50 time steps and a single reconstruction of the final state in TT format.

\textbf{Online efficiency over full time integration.} Hyper reduction reduceds the mean online time of TT-ROM for complete time integration from 2.536~s to 0.02180~s, leading to a $116.39\times$ speedup, while keeping the mean relative error at the same order of magnitude. 

\FloatBarrier
\subsection{Radiative Transfer Problem in a non-affine Lattice Medium}
\label{sec:5-5}
We consider the high-dimensional steady-state radiative transfer equation (RTE) with isotropic scattering~\cite{guo2026sweep}:
\begin{equation}
\label{eq:s5-d7}
\begin{gathered}
\Omega_x\partial_x\psi+\Omega_y\partial_y\psi+(\sigma_a+\sigma_s)\psi
=\sigma_s\rho+q,\\
\rho=\frac1{4\pi}\int_{S^2}\psi\,\mathrm d\Omega,
\end{gathered}
\end{equation}
Here, the unknown $\psi(x,y,\boldsymbol{\Omega};\boldsymbol{\mu})$ is the angular flux at the spatial position $(x,y)\in[0,5]^2$ in the unit direction $\boldsymbol{\Omega}=(\Omega_x,\Omega_y,\Omega_z)\in S^2$. The coefficients $\sigma_a(x,y;\boldsymbol{\mu})$ and $\sigma_s(x,y;\boldsymbol{\mu})$ are the absorption and scattering cross sections, respectively, and $q(x,y)$ is an isotropic source. The scalar flux $\rho(x,y;\boldsymbol{\mu})$ is the angularly averaged flux. The streaming term $\Omega_x\partial_x\psi+\Omega_y\partial_y\psi$ describes spatial transport. The term $(\sigma_a+\sigma_s)\psi$ accounts for absorption and scattering out of the direction $\boldsymbol{\Omega}$, while $\sigma_s\rho$ accounts for isotropic scattering into that direction. We impose zero inflow boundary conditions, corresponding to no incoming radiation.

\begin{figure}[htbp]
\centering
\includegraphics[width=\textwidth]{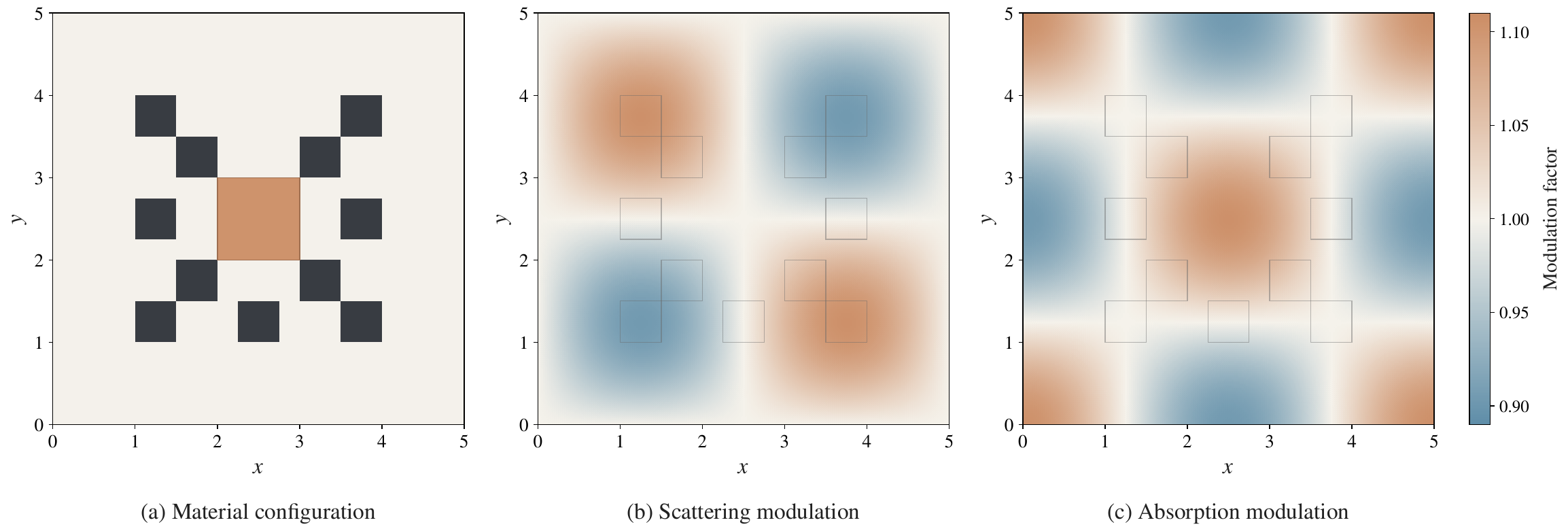}
\caption{Material configuration and spatial modulation factors for a radiative transfer problem in a lattice medium. The horizontal and vertical axes represent x and y, respectively. Subfigures (b) and (c) share a color scale centered at 1; the boundaries of the absorption region are marked by thin gray lines.}
\label{fig:3}
\end{figure}

For each $\boldsymbol{\mu}$, the anagular flux resides in a four-dimensional space--angle phase space (two spatial variables and two angles on $S^2$), motivating a space--angle low-rank representation.We vary $\boldsymbol{\mu}=(s,a,\alpha,\beta) \in [0.5,1.5] \times [95,105] \times [-0.2,0.2]^2$ around the nominal point $(1,100,0,0)$. Letting $\chi_B$ denote the absorption-region indicator, the cross sections are $\sigma_s=(1-\chi_B)s e^{\beta q_s}$ and $\sigma_a=\chi_B a e^{\alpha q_a}$, with $q_s(x,y)=\sin(2\pi x/5)\sin(2\pi y/5)$ and $q_a(x,y)=\cos(2\pi x/5)\cos(2\pi y/5)$. Figure~\ref{fig:3}(a) shows the $[0,5]^2$ layout (scattering background, dark-gray absorption squares, and central $(2,3)^2$ indicator source in orange), and Figs.~\ref{fig:3}(b)–(c) show $e^{\beta q_s}$ and $e^{\alpha q_a}$ for $\boldsymbol{\mu}^*=(1,100,0.1,-0.1)$.

A first-order tensor-product DG discretization uses \(60\times60\) cells, 288 angular directions, and five Gauss points per spatial direction for material integration. The full-order solver uses source iteration with diffusion synthetic acceleration (SI-DSA) \cite{adams2002fast}. It stops when the maximum absolute change in scalar flux is at most \(10^{-6}\). The training low-rank solver is the rank-adaptive sweep-based low-rank SI-DSA \cite{guo2026sweep} with the same outer-loop stopping criteria. The construction tolerance for the TT-reduced basis is \(10^{-5}\), while the construction tolerance and local rank threshold for the TT-DEIM basis for each material field are both \(10^{-3}\). The resulting solution space has dimensions ranging from 27 to 30 and both material DEIM spaces have dimensions of 3-4.

\textbf{Treatment of material terms.} TT-cross-DEIM~\cite{ghahremani2024cross}, with the rank criterion described in Appendix~\ref{app:auxiliary-cross}, approximates the smooth factors \(se^{(\beta q_s)}\) and \(ae^{(\alpha q_a)}\), excluding the fixed geometric masks. Combining the resulting TT-DEIM bases with the masks and DG integrals gives operator components whose reduced projections are precomputed. The solution error is measured using the relative Euclidean error defined in Section~\ref{sec:5-1}, applied to the vectorized space--angle coefficients, without applying the DG mass matrix or angular quadrature weights.

\textbf{Partial entry access.} Table~\ref{tab:rte-entry-access} shows that TT-cross-DEIM matches the accuracy order of full-access TT-SVD (at relative tolerance $10^{-3}$) while sampling fewer than 10\% of tensor entries: 9.69\% for scattering (relative Frobenius error $3.206\times10^{-4}$) and 9.87\% for absorption ($2.800\times10^{-4}$), with a maximum TT rank of 5 versus 4 for TT-SVD. Thus, smooth material fields can be accurately compressed from partial evaluations when the material geometry is fixed.

\begin{table}[htbp]
\centering
\caption{Comparison of Entry Access and Approximation Accuracy for Radiative Transfer Material Tensors}
\label{tab:rte-entry-access}
\footnotesize
\setlength{\tabcolsep}{4pt}
\renewcommand{\arraystretch}{1.25}
\begin{adjustbox}{max width=\linewidth}
\begin{tabular}{@{}llccc@{}}
\toprule
Material & Construction method & Entry access(\%) & \shortstack{Mean rel. Frobenius error} & \shortstack{Max TT rank}\\
\midrule
\multirow{2}{*}{Scattering $g_s$} & TT-cross-DEIM & 9.69 & $3.206\times10^{-4}$ & 5\\
& Full-access TT-SVD & 100.00 & $2.969\times10^{-4}$ & 4\\
\midrule
\multirow{2}{*}{Absorption $g_a$} & TT-cross-DEIM & 9.87 & $2.800\times10^{-4}$ & 5\\
& Full-access TT-SVD & 100.00 & $2.994\times10^{-4}$ & 4\\
\bottomrule
\end{tabular}
\end{adjustbox}
\par\smallskip\begin{minipage}{\textwidth}\footnotesize\raggedright \end{minipage}
\end{table}

\textbf{Online efficiency.} Hyper reduction reduces the mean online time from 9.273~s (TT-Galerkin ROM) to $6.037\times10^{-3}$~s, a $1532.37\times$ mean speedup, while maintaining comparable mean relative errors ($4.762\times10^{-5}$ vs.\ $1.560\times10^{-5}$). Unlike the Galerkin ROM, which assembles and projects the parameter-dependent spatial material matrix at every test parameter, the hyper-reduced TT-ROM recovers expansion coefficients from sampled material values and combines precomputed reduced matrices. 

\FloatBarrier
\Needspace{6\baselineskip}
\subsection{Summary of Accuracy, Storage, and Computational Cost}
\label{sec:5-6}

We now compare the four examples using the protocol in Section~\ref{sec:5-1}. The low-rank solver provides the baseline for assessing the additional online savings from model reduction, while the comparison between the two ROMs measures the benefit of hyper reduction.

\textbf{Low-rank solutions versus full-order references.} Table~\ref{tab:lr-solutions} lists the accuracy and storage characteristics of low-rank solutions on the test set. The mean relative error ranges from \(6.718\times10^{-8}\) to \(5.111\times10^{-5}\), and the average storage compression ratio ranges from 1.74 to 5.37, showing that the low-rank solvers produce accurate solutions in compressed form. These ratios measure solution storage only. 

\begin{table}[htbp]
\centering
\caption{Accuracy and Solution Storage of Low-Rank Solutions on the Test Set}
\label{tab:lr-solutions}
\footnotesize
\setlength{\tabcolsep}{4pt}
\renewcommand{\arraystretch}{1.25}
\begin{adjustbox}{max width=\linewidth}
\begin{tabular}{@{}lllcccc@{}}
\toprule
Problem & FOM solver & LR solver & \shortstack{Mean rel.\\error} & \shortstack{Max rel.\\error} & \shortstack{Mean TT\\rank} & \shortstack{Mean storage\\compression}\\
\midrule
Parametric diffusion & Sparse LU & Preconditioned lrAA & $4.831\times10^{-6}$ & $2.637\times10^{-5}$ & 9.96 & $3.14\times$\\
Parametric Bratu & Sparse Newton--LU & Cross-DEIM lrAA & $6.718\times10^{-8}$ & $5.683\times10^{-7}$ & 18.02 & $1.74\times$\\
Allen--Cahn & Newton--GMRES & Preconditioned lrAA & $5.111\times10^{-5}$ & $1.654\times10^{-4}$ & 26.73 & $5.37\times$\\
Lattice-medium RTE & SI-DSA & Low-rank SI-DSA & $3.364\times10^{-6}$ & $6.284\times10^{-6}$ & 138.93 & $2.03\times$\\
\bottomrule
\end{tabular}
\end{adjustbox}
\par\smallskip\begin{minipage}{\textwidth}\footnotesize\raggedright Note:  SI denotes source iteration; DSA denotes diffusion synthetic acceleration. \end{minipage}
\end{table}

\textbf{ROM accuracy and online efficiency.} Table~\ref{tab:rom-performance} presents the key results of this section: hyper reduction significantly reduces online time for all four problem types, with mean paired speedups of 8.42 (diffusion), 10.39 (Bratu), 116.39 (Allen--Cahn), and 1,532.37 (radiative transfer), respectively. The hyper-reduced TT-ROM maintains mean relative errors comparable to those of the TT Galerkin ROM, on the order of \(10^{-5}\)--\(10^{-4}\) for all four cases, with maximum errors below \(5\times10^{-3}\). Both TT-ROMs also reduce computational time relative to the low-rank solver at the reported error levels. These timings exclude offline construction; the total-cost comparison below accounts for the additional expense of hyper reduction.

\begin{table}[htbp]
\centering
\caption{Accuracy and Online Efficiency of the TT Galerkin ROM and the TT Hyper-ROM}
\label{tab:rom-performance}
\footnotesize
\setlength{\tabcolsep}{2.5pt}
\renewcommand{\arraystretch}{1.2}
\medmuskip=2mu\relax\thickmuskip=3mu\relax
\begin{adjustbox}{max width=\textwidth}
\begin{tabular}{@{}llccccc@{}}
\toprule
Problem & Method & \shortstack{Online time (s)\\mean $\pm$ std} & \shortstack{Mean rel. error\\mean $\pm$ std} & \shortstack{Max rel.\\error} & \shortstack{Speedup\\vs LR} & \shortstack{Speedup\\vs Galerkin}\\
\midrule
\multirow{3}{*}{Parametric diffusion}
& LR solver & $(2.760\pm0.000)\times10^{-1}$ & $(4.831\pm0.000)\times10^{-6}$ & $2.637\times10^{-5}$ & $1.00\times$ & ---\\
& TT Galerkin & $(3.359\pm0.048)\times10^{-3}$ & $(1.676\pm0.160)\times10^{-4}$ & $1.865\times10^{-3}$ & $82.37\times$ & $1.00\times$\\
& TT hyper-ROM & $(4.118\pm0.102)\times10^{-4}$ & $(2.280\pm0.273)\times10^{-4}$ & $4.854\times10^{-3}$ & $691.80\times$ & $8.42\times$\\
\midrule
\multirow{3}{*}{Parametric Bratu}
& LR solver & $(7.845\pm0.000)\times10^{-2}$ & $(6.718\pm0.000)\times10^{-8}$ & $5.683\times10^{-7}$ & $1.00\times$ & ---\\
& TT Galerkin & $(1.697\pm0.091)\times10^{-3}$ & $(1.905\pm0.385)\times10^{-4}$ & $1.156\times10^{-3}$ & $46.57\times$ & $1.00\times$\\
& TT hyper-ROM & $(1.630\pm0.023)\times10^{-4}$ & $(2.825\pm0.803)\times10^{-4}$ & $3.144\times10^{-3}$ & $481.61\times$ & $10.39\times$\\
\midrule
\multirow{3}{*}{Allen--Cahn}
& LR solver & $(1.839\pm0.000)\times10^{1}$ & $(5.111\pm0.000)\times10^{-5}$ & $1.654\times10^{-4}$ & $1.00\times$ & ---\\
& TT Galerkin & $(2.536\pm0.036)\times10^{0}$ & $(4.314\pm0.277)\times10^{-5}$ & $1.586\times10^{-4}$ & $7.26\times$ & $1.00\times$\\
& TT hyper-ROM & $(2.180\pm0.028)\times10^{-2}$ & $(8.326\pm1.670)\times10^{-5}$ & $6.369\times10^{-4}$ & $843.84\times$ & $116.39\times$\\
\midrule
\multirow{3}{*}{Lattice-medium RTE}
& LR solver & $(2.306\pm0.000)\times10^{2}$ & $(3.364\pm0.000)\times10^{-6}$ & $6.284\times10^{-6}$ & $1.00\times$ & ---\\
& TT Galerkin & $(9.273\pm1.231)\times10^{0}$ & $(1.560\pm0.309)\times10^{-5}$ & $6.748\times10^{-5}$ & $25.22\times$ & $1.00\times$\\
& TT hyper-ROM & $(6.037\pm0.471)\times10^{-3}$ & $(4.762\pm0.664)\times10^{-5}$ & $1.520\times10^{-4}$ & $38421.58\times$ & $1532.37\times$\\
\bottomrule
\end{tabular}
\end{adjustbox}
\par\smallskip\begin{minipage}{\textwidth}\footnotesize\raggedright Note: TT Galerkin denotes the TT Galerkin ROM. Both speedup columns report mean paired speedups. ``Mean \(\pm\) standard deviation'' is calculated from the test means of each of the five training seeds, using the sample standard deviation.\end{minipage}
\end{table}

\textbf{Offline cost and amortization.} Hyper reduction adds offline costs for TT approximation, DEIM basis construction, sampling-point selection, and projection. Shared costs of snapshot generation, solution-basis construction, and constant-term projection cancel in the comparison with TT Galerkin. For $M\in\{\mathrm G,\mathrm H\}$, the cost of $Q$ queries and the break-even query count are

\begin{equation}
\label{eq:s5-d9}
\begin{aligned}
C_M(Q)&=T_{\mathrm{off}}^M+Q\overline T_{\mathrm{on}}^M,\qquad N_{\mathrm{BE}}=\left\lceil
\frac{T_{\mathrm{off}}^{\mathrm H}-T_{\mathrm{off}}^{\mathrm G}}
{\overline T_{\mathrm{on}}^{\mathrm G}-\overline T_{\mathrm{on}}^{\mathrm H}}
\right\rceil.
\end{aligned}
\end{equation}

Table~\ref{tab:offline-cost} shows that the additional hyper reduction cost is recovered after tens of queries for Allen--Cahn and RTE, while Bratu and diffusion require hundreds to thousands of queries.

\begin{table}[htbp]
\centering
\caption{Additional Offline Cost and Break-Even Queries Relative to TT Galerkin}
\label{tab:offline-cost}
\footnotesize
\setlength{\tabcolsep}{4pt}
\renewcommand{\arraystretch}{1.25}
\begin{adjustbox}{max width=\linewidth}
\begin{tabular}{@{}lccc@{}}
\toprule
Problem & Additional offline cost (s)
        & Online saving per query (s)
        & Break-even queries \\
\midrule
Parametric diffusion & 21.89474 & $2.94744\times10^{-3}$ & 7429\\
Parametric Bratu & 1.25561 & $1.53374\times10^{-3}$ & 819\\
Allen--Cahn & 162.26930 & 2.51445 & 65\\
Lattice-medium RTE & 91.31278 & 9.26657 & 10\\
\bottomrule
\end{tabular}
\end{adjustbox}
\par\smallskip\begin{minipage}{\textwidth}\footnotesize\raggedright \end{minipage}
\end{table}

\FloatBarrier

\section{Conclusions}
\label{sec:6}
To combine the efficiency of low-rank PDE solvers with data-driven ROM for nonlinear and nonaffine parametric PDEs, we developed a hyper-reduced TT-ROM constructed directly from streaming TT data. The framework combines incremental low-rank TT compression with nested TT-DEIM to build compressed reduced bases and select nested interpolation points directly from streaming TT data through only core-wise operations. Online, evaluations of the nonaffine and nonlinear terms at these fixed points are combined with precomputed reduced projections to assemble and solve the reduced system. 

Numerical experiments on nonaffine diffusion, the Bratu problem, the Allen--Cahn equation, and a non-affine radiative transport demonstrated that the hyper-reduced models can maintain accuracy comparable to TT-Galerkin models without hyper-reduction while substantially reducing online computational time.

\appendix

\section{TT-Cross Approximation of Training Tensors}
\label{app:auxiliary-cross}
 We briefly review the iterative TT-CUR-DEIM in \cite{ghahremani2024cross}, including alternating DEIM sweeps and GappyPOD+E enrichment~\cite{peherstorfer2020stability}. We specify below a modified local rank criterion for fixed training tensors and its associated sampling and termination rules.

\textbf{Sampling and reconstruction.}  Let $\mathcal{X}\in\mathbb R^{q_1\times\cdots\times q_d}$ be the target tensor, with entries evaluated on demand. At ranks $(1,r_1,\ldots,r_{d-1},1)$, nested index sets satisfy $|\mathcal{I}_k|=|\mathcal{J}_k|=r_k$ and
\begin{equation}
\label{eq:s3-d9}
\mathcal{I}_k\subseteq\mathcal{I}_{k-1}\times[q_k],\qquad
\mathcal{J}_k\subseteq[q_{k+1}]\times\mathcal{J}_{k+1},\qquad
\mathcal{I}_0=\mathcal{J}_d=\{()\}.
\end{equation}
At fixed ranks, alternate DEIM sweeps until the ordered right index sets stabilize or the sweep limit is reached, then reconstruct the TT cores in a final left-to-right sweep with fixed suffixes~\cite{ghahremani2024cross}.

\textbf{Local rank criterion.} Retain the existing suffixes and enrich them from right to left, starting with $\mathcal{J}_d^+=\{()\}$. At interface $k$, the candidates are $[q_{k+1}]\times\mathcal{J}_{k+1}^+$. Apply the GappyPOD+E \cite{peherstorfer2020stability} score to the $r_k$-dimensional column basis of
\[
\operatorname{unfold}_R\bigl(\mathcal{X}(\mathcal{I}_k,:,\mathcal{J}_{k+1}^+)\bigr)^T,
\]
using the rows corresponding to $\mathcal{J}_k$ as the initial selection. Add up to $m_{\mathrm{os}}\geq1$ unused sampling multi-indices to obtain $\mathcal{J}_k^+\supseteq\mathcal{J}_k$. For $r_k=1$, select additional rows by the magnitude of the basis-vector entries.

Evaluate the enlarged local matrix
\begin{equation}
\label{eq:s3-d14}
\boldsymbol{B}_k^+=\operatorname{unfold}_L\bigl(\mathcal{X}(\mathcal{I}_{k-1},:,\mathcal{J}_k^+)\bigr)
\in\mathbb R^{(r_{k-1}q_k)\times|\mathcal J_k^+|}.
\end{equation}
The implementation uses the first singular value omitted by rank $r_k$:
\begin{equation}
\label{eq:s3-d15}
\epsilon_k=
\begin{cases}
\displaystyle\frac{\sigma_{r_k+1}(\boldsymbol{B}_k^+)}{\|\boldsymbol{B}_k^+\|_F},
& r_k<\min\{r_{k-1}q_k,|\mathcal{J}_k^+|\},\quad\|\boldsymbol{B}_k^+\|_F>0,\\[2pt]
0,&\text{otherwise}.
\end{cases}
\end{equation}
The criterion in~\cite[Section~2.6]{ghahremani2024cross} uses the smallest retained singular value to guide rank increases and decreases during time integration. Our test detects an additional sampled direction in a fixed tensor and flags rank $r_k$ for an increase when $\epsilon_k>\varepsilon_{\mathrm{rank}}$.

\textbf{Rank updates and termination.} Increase each detected rank by at most one, subject to
\[
\widetilde r_k\in\{r_k,r_k+1\},\qquad
1\leq\widetilde r_k\leq
\min\{r_{\max},\widetilde r_{k-1}q_k,q_{k+1}\widetilde r_{k+1}\},
\qquad \widetilde r_0=\widetilde r_d=1.
\]
Initialize the next reconstruction by retaining valid old suffixes, adding available enriched suffixes, and filling any remaining positions from the nested candidates. Repeat reconstruction, enrichment, and testing until all indicators satisfy the threshold or a rank or iteration limit prevents further adaptation. The settings $\Theta_{\mathrm{approx}}$ comprise $\varepsilon_{\mathrm{rank}}$, the initial ranks, $m_{\mathrm{os}}$, the sweep and adaptation limits, and $r_{\max}$.

The reconstruction assumes full-column-rank sampled matrices and nonsingular cross matrices. The local rank test can miss unsampled directions and does not certify a full-tensor error bound; termination at a limit may also leave the local threshold unsatisfied. The analysis in Section~\ref{sec:4} assumes a bound on the supplied TT snapshot errors and applies independently of this approximation routine.

\section*{Declaration of Generative AI Use}
During the preparation of this manuscript, the authors used ChatGPT and Gemini to assist with code implementation and language polishing. The authors take full responsibility for the content of the manuscript.

\bibliographystyle{siam}
\bibliography{ref}

\end{document}